\documentclass{article}
\usepackage{amsmath,amssymb}
\makeindex

\begin{document}

\title{A Lefschetz formula for higher rank}
\author{Anton Deitmar}

\date{}
\maketitle

\setlength{\parskip}{7pt}

\def \1{{\bf 1}}
\def \a{{{\mathfrak a}}}
\def \ad{{\rm ad}}
\def \al{\alpha}
\def \ar{{\alpha_r}}
\def \A{{\mathbb A}}
\def \Ad{{\rm Ad}}
\def \Aut{{\rm Aut}}
\def \b{{{\mathfrak b}}}
\def \bs{\backslash}
\def \B{{\cal B}}
\def \c{{\mathfrak c}}
\def \cent{{\rm cent}}
\def \conj{{\rm conj}}
\def \C{{\mathbb C}}
\def \CA{{\cal A}}
\def \CB{{\cal B}}
\def \CC{{\cal C}}
\def \CE{{\cal E}}
\def \CF{{\cal F}}
\def \CG{{\cal G}}
\def \CH{{\cal H}}
\def \CHC{{\cal HC}}
\def \CL{{\cal L}}
\def \CM{{\cal M}}
\def \CN{{\cal N}}
\def \CP{{\cal P}}
\def \CQ{{\cal Q}}
\def \CO{{\cal O}}
\def \CS{{\cal S}}
\def \CT{{\cal T}}
\def \CV{{\cal V}}
\def \det{{\rm det}}
\def \df{\ \begin{array}{c} _{\rm def}\\ ^{\displaystyle =}\end{array}\ }
\def \diag{{\rm diag}}
\def \dist{{\rm dist}}
\def \End{{\rm End}}
\def \eps{\varepsilon}
\def \Fx{{\mathfrak x}}
\def \FX{{\mathfrak X}}
\def \g{{{\mathfrak g}}}
\def \ga{\gamma}
\def \G{{\mathbb G}}
\def \Ga{\Gamma}
\def \GL{{\rm GL}}
\def \h{{{\mathfrak h}}}
\def \Hom{{\rm Hom}}
\def \im{{\rm im}}
\def \Im{{\rm Im}}
\def \Ind{{\rm Ind}}
\def \k{{{\mathfrak k}}}
\def \K{{\cal K}}
\def \l{{\mathfrak l}}
\def \la{\lambda}
\def \lap{\triangle}
\def \La{\Lambda}
\def \m{{{\mathfrak m}}}
\def \mod{{\rm mod}}
\def \n{{{\mathfrak n}}}
\def \name{\bf}
\def \N{\mathbb N}
\def \o{{\mathfrak o}}
\def \ord{{\rm ord}}
\def \O{{\cal O}}
\def \p{{{\mathfrak p}}}
\def \ph{\varphi}
\def \Pr{{\rm Pr}}
\def \prf{\noindent{\bf Proof: }}
\def \Per{{\rm Per}}
\def \q{{\mathfrak q}}
\def \qed{\ifmmode\eqno $\square$\else\noproof\vskip 12pt plus 3pt minus 9pt \fi}
 \def\noproof{{\unskip\nobreak\hfill\penalty50\hskip2em\hbox{}%
     \nobreak\hfill $\square$\parfillskip=0pt%
     \finalhyphendemerits=0\par}}
\def \Q{\mathbb Q}
\def \res{{\rm res}}
\def \R{{\mathbb R}}
\def \Re{{\rm Re \hspace{1pt}}}
\def \r{{\mathfrak r}}
\def \ra{\rightarrow}
\def \rank{{\rm rank}}
\def \Rep{{\rm Rep}}
\def \setminus{\begin{picture}(18,10)\put(4,6)
                {\line(2,-1){10}}\end{picture}}
\def \sm{{\rm sm}}
\def \supp{{\rm supp}}
\def \Spin{{\rm Spin}}
\def \t{{{\mathfrak t}}}
\def \T{{\mathbb T}}
\def \tr{{\hspace{1pt}\rm tr\hspace{2pt}}}
\def \vol{{\rm vol}}
\def \z{\zeta}
\def \Z{\mathbb Z}
\def \={\ =\ }
\def \({\left(}
\def \){\right)}

\newcommand{\choice}[4]{\left\{
            \begin{array}{cl}#1&#2\\ #3&#4\end{array}\right.}
\newcommand{\rez}[1]{\frac{1}{#1}}
\newcommand{\der}[1]{\frac{\partial}{\partial #1}}
\newcommand{\norm}[1]{\parallel #1 \parallel}
\renewcommand{\matrix}[4]{\left( \begin{array}{cc}#1 & #2 \\ #3 & #4 \end{array}
            \right)}

\newtheorem{theorem}{Theorem}[section]
\newtheorem{conjecture}[theorem]{Conjecture}
\newtheorem{lemma}[theorem]{Lemma}
\newtheorem{corollary}[theorem]{Corollary}
\newtheorem{proposition}[theorem]{Proposition}

$$ $$

\section*{Abstract}
In this paper a Lefschetz formula is proved for the geodesic flow of a compact locally symmetric space.
The flow is described in terms of actions of split tori of various dimensions and the geometric side of the Lefschetz formula is a sum over closed geodesics which correspond to a given torus.
The cohomological side is given in terms of Lie algebra cohomology.

\tableofcontents

\section*{Introduction}

In this paper we prove a Lefschetz formula for the geodesic flow of a compact locally symmetric space.
The first case of such a formula appears in \cite{Juhl}, were it is proved for compact quotients of symmetric spaces of rank one.
In higher rank, the geodesic flow extends to an action of a higher dimensional torus on the sphere bundle. The sphere bundle itself carries a stratification according to the orbit type and the Lefschetz formula is formulated for each stratum separately.

In the first section we fix notation and collect some prerequisites, most of them going back to Harish-Chandra's work.
In section 2 we construct Euler-Poincar\'e functions in a more general setting as in the literature. The Selberg trace formula is given in section 3.
Finally, in section 4 we formulate and prove the Lefschetz formula which connects geometric information on closed orbits with spectral data from the action of the flow on Lie algebra cohomology.

This formula can be used to show that Zeta functions of Selberg type have an analytic continuation if the torus is one dimensional (within a higher rank group).
But the Lefschetz formula gives valuable information for each dimension.
A weaker version of the highest dimensional case has been used in \cite{primgeo} to derive asymptotic formulae for the length distribution of closed geodesics.

\newpage
\section{Prerequisites}

In this section we collect some facts from literature which
will be needed in the sequel. Proofs will only be given by
sketches or references.

\subsection{Notations} \label{notations}

We denote Lie groups by upper case roman letters $G,H,K,\dots$ and
the corresponding real Lie algebras by lower case German letters
with index $0$, that is: $\g_0,\h_0,\k_0,\dots$. The complexified
Lie algebras will be denoted by $\g,\h,\k,\dots$, so, for
example: $\g =\g_0\otimes_\R \C$.

For a Lie group $L$ with Lie algebra $\l_0$ let $\Ad : L\ra
GL(\l_0)$ be the adjoint representation (\cite{bourb} III.3.12).
By definition, $\Ad(g)$ is the differential of the map $x\mapsto
gxg^{-1}$ at the point $x=e$. Then $\Ad(g)$ is a Lie algebra
automorphism of $\g_0$. A Lie group $L$ is said to be {\it of
inner type} if $\Ad(L)$ lies in the complex
adjoint group of the Lie algebra $\l$.

A real Lie group $G$ is said to be a {\it real reductive
group} if there is a linear algebraic
group $\CG$ defined over $\R$ which is reductive as an algebraic
group \cite{Borel-lingroups} and a morphism $\alpha : G\ra \CG(\R)$
with finite kernel and cokernel (Since we do not insist that the
image of $\alpha$ is normal the latter condition means that ${\rm
im}(\alpha)$ has finite index in $\CG(\R)$). This implies in
particular that $G$ has only finitely many connected components
(\cite{Borel-lingroups} 24.6.c).

A real reductive group $G$ has a maximal compact subgroup $K$ which meets every connected component.
The group $G$ is of inner type if and only if $\Ad(K)$ lies in the complex adjoint group of $\g$.

Note that any real reductive group $G$ of inner type is of {\it
Harish-Chandra class}, i.e., $G$ is of
inner type, the Lie algebra $\g_0$ of $G$ is reductive, $G$ has
finitely many connected components and the connected subgroup
$G_{der}^0$ corresponding to the Lie subalgebra
$[\g_0,\g_0]$ has finite center.

The following are of importance:
\begin{itemize}
\item
a connected semisimple Lie group with finite center is a real
reductive group of inner type (\cite{wall-rg1} 2.1.3),
\item
if $G$ is a real reductive group of inner type and $P=MAN$ is the
Langlands decomposition (\cite{wall-rg1} 2.2.7) of a parabolic
subgroup then the groups $M$ and $AM$ are real reductive of inner
type (\cite{wall-rg1} 2.2.8)
\end{itemize}

The usual terminology of algebraic groups carries over to real
reductive groups, for example a {\it torus} (or a
Cartan subgroup) in $G$ is the inverse image of (the real points
of) a torus or a Cartan subgroup in $\CG(\R)$. An element of $G$
will be called {\it semisimple} if it
lies in some torus of $G$. The {\it split component} of $G$ is the
identity component of the greatest split torus in the center of
$G$. Note that for a real reductive group the Cartan subgroups are
precisely the maximal tori.

Let $G$ be a real reductive group then there exists a {\it Cartan
involution} i.e., an automorphism
$\theta$ of $G$ satisfying $\theta^2=Id$ whose
fixed point set is a maximal compact subgroup $K$ and which is the
inverse ($a\mapsto a^{-1}$) on the split component of $G$. All
Cartan involutions are conjugate under automorphisms of $G$.

Fix a Cartan involution $\theta$ with fixed point set equal to the maximal
compact subgroup $K$ and let $\k_0$ be the Lie algebra of $K$.
 The group $K$ acts on $\g_0$ via the adjoint
representation and there is a $K$-stable decomposition
$\g_0=\k_0\oplus\p_0$, where $\p_0$ is the
eigenspace of (the differential of) $\theta$ to the eigenvalue
$-1$. Write $\g=\k\oplus\p$ for the complexification. This is
called the {\it Cartan decomposition}.

\begin{lemma} There is a symmetric bilinear form $B : \g_0\times\g_0\ra \R$
such that
\begin{itemize}
\item
$B$ is invariant, that is $B(\Ad(g)X,\Ad(g)Y)=B(X,Y)$ for all
$g\in G$ and all $X,Y\in\g_0$ and
\item
$B$ is negative definite on $\k_0$ and positive definite on its
orthocomplement $\p_0=\k_0^\perp\subset\g_0$.
\end{itemize}
\end{lemma}

\prf For $X\in\g_0$ let $\ad(X):\g_0\ra\g_0$ be the adjoint
defined by $\ad(X)Y=[X,Y]$. Since for $g\in G$ the map $\Ad(g)$ is
a Lie algebra homomorphism we infer that
$\ad(\Ad(g)X)=\Ad(g)\ad(X)\Ad(g)^{-1}$ and therefore the {\it
Killing form}
$$
B_K(X,Y)\= \tr (\ad(X)\ad(Y))
$$
is invariant. It is known that if $B_K$ is nondegenerate, i.e.,
$\g_0$ is semisimple, then $B=B_K$ satisfies the claims of the
lemma. In the general case we have
$\g_0=\a_0\oplus\c_0\oplus\g_0'$, where $\a_0\oplus\c_0$ is the
center of $\g_0$ and $\g_0'$ its derived algebra, which is
semisimple and so $B_K|_{\g_0'}$ is nondegenerate, whereas
$B_K|_{\a_0\oplus\c_0}=0$. Further $\c_0$ is the eigenspace of
$\theta$ in the center of $\g_0$ corresponding to the eigenvalue
$1$ whereas $\a_0$ is the eigenspace of $-1$. For $g\in G$ the
adjoint $\Ad(g)$ is easily seen to preserve $\a_0$ and $\c_0$, so
we get a representation $\rho:G\ra GL(\a_0)\times GL(\c_0)$. This
representation is trivial on the connected component $G^0$ of $G$
hence it factors over the finite group $G/G^0$. Therefore there is
a positive definite symmetric bilinear form $B_a$ on $\a_0$ which
is invariant under $\Ad_c$ and similarly a negative definite
symmetric bilinear form $B_c$ on $\c_0$ which is invariant. Let
$$
B\= B_a\oplus B_c\oplus B_K|_{\g_0'}.
$$
Then $B$ satisfies the claims of the lemma. \qed

Let $U(\g_0)$ be the universal enveloping algebra
of $\g_0$. It can be constructed as the quotient of the tensorial
algebra
$$
T(\g_0) \= \R\oplus\g_0\oplus (\g_0\otimes\g_0)\oplus\dots
$$
by the two-sided ideal generated by all elements of the form
$X\otimes Y-Y\otimes X-[X,Y]$, where $X,Y\in\g_0$.

The algebra $U(\g_0)$ can be identified with the $\R$-algebra of
all left invariant differential operators on $G$. Let $\g
=\g_0\otimes\C$ be the complexification of $\g_0$ and
$U(\g)$ be its universal enveloping algebra which is the same as
the complexification of $U(\g_0)$. Then $\g$ is a subspace of
$U(\g)$ which generates $U(\g)$ as an algebra and any Lie algebra
representation of $\g$ extends uniquely to a representation of
the associative algebra $U(\g)$.

The form $B$ gives an identification of $\g_0$ with its dual space
$\g_0^*$. On the other hand $B$ defines an element in
$\g_0^*\otimes \g_0^*$. Thus we get a natural element in
$\g_0\otimes\g_0\subset T(\g_0)$. The image $C$ of this element
in $U(\g_0)$ is called the {\it Casimir operator} attached to $B$.
 It is a differential operator of order
two and it lies in the center of $U(\g_0)$. In an more concrete
way the Casimir operator can be described as follows: Let $X_1,
\dots, X_m$ be a basis of $\g_0$ and let $Y_1, \dots, Y_m$ be the
dual basis with respect to the form $B$  then the Casimir operator
is given by
$$
C\= X_1 Y_1 +\dots +X_m Y_m.
$$
We have

\begin{lemma}\label{C-central}
For any $g\in G$ the Casimir operator is invariant under
$\Ad(g)$, that is $\Ad(g) C=C$.
\end{lemma}

\prf This follows from the invariance of $B$
\qed

Let $X$ denote the quotient manifold $G/K$. The tangent space at
$eK$ identifies with $\p_0$ and the form $B$ gives a
$K$-invariant positive definite inner product on this space.
Translating this by elements of $G$ defines a $G$-invariant
Riemannian metric on $X$. This makes $X$ the most general
globally symmetric space of the noncompact type \cite{helg}.

Let $\hat{G}$ denote the {\it unitary
dual} of $G$, i.e., $\hat{G}$ is the set of
isomorphism classes of irreducible unitary representations of $G$.

Let $(\pi,V_\pi)$ be a continuous representation of $G$ on some
Banach space $V_\pi$. The subspace $V_\pi^\infty$ of {\it smooth
vectors} is defined to be the subspace of
$V_\pi$ consisting of all $v\in V_\pi$ such that the map
$g\mapsto \pi(g)v$ is smooth. The universal enveloping algebra
$U(\g)$ operates on $V_\pi^\infty$ via
$$
\pi(X): v\mapsto X_g(\pi(g)v)\ |_{g=e}
$$
for $X$ in $\g$.

A {\it $(\g,K)$-module} is by definition a
complex vector space $V$ which is a $K$-module such that for each
$v\in V$ the space spanned by the orbit $K.v$ is finite
dimensional. Further $V$ is supposed to be a $\g$-module and the
following compatibility conditions should be satisfied:
\begin{itemize}
\item
for $Y\in\k\subset\g$ and $v\in V$ it holds
$$
Y.v\= \left.\frac{d}{dt}\right|_{t=0} \exp(tY).v,
$$
\item
for $k\in K$, $X\in\g$ and $v\in V$ we have
$$
k.X.v\= \Ad(k)X.k.v.
$$
\end{itemize}
A $(\g,K)$-module $V$ is called {\it irreducible} or simple if it
has no proper submodules and it is called of finite length if
there is a finite filtration
$$
0=V_0\subset\dots\subset V_n=V
$$
of submodules such that each quotient $V_j/V_{j-1}$ is
irreducible. Further $V$ is called {\it admissible} if for each
$\tau\in\hat{K}$ the space $\Hom_K(\tau,V)$ is finite dimensional.
An admissible $(\g,K)$-module of finite length is called a {\it
Harish-Chandra module}.

Now let $(\pi,V_\pi)$ again be a Banach representation of $G$. 
For each $\tau\in\hat K$ let $V_\pi(\tau)$ denote the isotypical component of $\tau$, i.e., $V_\pi(\tau)$ is the image of the map
\begin{eqnarray*}
\Hom_K(V_\tau,V_\pi)\otimes V_\tau &\to& V_\pi,\\
(\alpha,v) &\mapsto& \alpha(v).
\end{eqnarray*}
Let
$V_{\pi,K}$ be the subspace of $V_\pi$ consisting of all vectors
$v\in V_\pi$ such that the $K$-orbit $\pi(K)v$ spans a finite
dimensional space. Then $V_{\pi,K}$  is called the space of {\it
$K$-finite vectors} in $V_\pi$. The
space $V_{\pi,K}$ is no longer a $G$-module but remains a
$K$-module. Further the space $V_{\pi,K}^\infty = V_{\pi,K}\cap
V_\pi^\infty$ is dense in $V_{\pi,K}$ and is stable under $K$, so
$V_{\pi,K}^\infty$ is a $(\g,K)$-module. By abuse of notation we
will often write $\pi$ instead of $V_\pi$ and $\pi_K^\infty$
instead of $V_{\pi,K}^\infty$. The representation $\pi$ is called
admissible if $\pi_K^\infty$ is. In that case we have
$\pi_K^\infty =\pi_K$ since a dense subspace of a finite
dimensional space equals the entire space. So then $\pi_K$ is a
$(\g,K)$-module.
Further, the representation $\pi$ is called a \emph{Harish-Chandra representation} if $\pi_K$ is a Harish-Chandra module.

\begin{lemma}
Let $(\pi,V_\pi)$ be an irreducible admissible representation of
$G$ then the Casimir operator $C$ acts on $V_\pi^\infty$ by a
scalar denoted $\pi(C)$.
\end{lemma}

\prf By the formula $\pi(g)\pi(C)\pi(g)^{-1}=\pi(\Ad(g)C)$ and
Lemma \ref{C-central} we infer that $\pi(C)$ commutes with
$\pi(g)$ for every $g\in G$. Therefore the claim follows from the
Lemma of Schur (\cite{wall-rg1} 3.3.2).  \qed

Any $f\in L^1(G)$ will define a continuous operator
$$
\pi(f)\=\int_Gf(x)\pi(x) dx
$$
on $V_\pi$.

Let $N$ be a natural number and let
$L_{2N}^1(G)$ be the set of all $f\in
C^{2N}(G)$ which satisfy $Df\in L^1(G)$ for any $D\in U(\g)$ with
$\deg(D)\le 2N$.

\begin{lemma}
Let $N$ be an integer $>\frac{\dim G}{2}$. Let $f\in L_{2N}^1(G)$
then for any irreducible unitary representation $\pi$ of $G$ the
operator $\pi(f)$ will be of trace class.
\end{lemma}

\prf Let $C$ denote the Casimir operator of $G$ and let $C_K$ be
the Casimir operator of $K$. Let $\lap = -C+2C_K\in U(\g)$ the
group Laplacian. It is known that for some $a>0$ the operator
$\pi(\lap +a)$ is positive and $\pi(\lap +a)^{-N}$ is of trace
class. Let $g=(\lap +a)^Nf$ then $g\in L^1(G)$, so $\pi(g)$ is
defined and gives a continuous linear operator on $V_\pi$. We
infer that $\pi(f)=\pi(\lap +a)^{-N}\pi(g)$ is of trace class.
\qed

Finally we need some more notation. The form $\langle X,Y\rangle
=-B(X,\theta(Y))$ is positive definite on $\g_0$ and therefore
induces a positive definite left invariant top differential form
$\omega_L$ on any closed subgroup $L$ of $G$. If $L$ is compact we
set
$$
v(L)\= \int_L\omega_L.
$$
Let $H=AB$ be a $\theta$-stable Cartan subgroup where $A$ is the
connected split component of $H$ and $B$ is compact. The double
use of the letter $B$ here will not cause any confusion. Then
$B\subset K$. Let $\Phi$ denote the root system of $(\g,\h)$,
where $\g$ and $\h$ are the complexified Lie algebras of $G$ and
$H$. Let $\g =\h\oplus\bigoplus_{\alpha\in\Phi}\g_\alpha$ be the
root space decomposition. Let $x\ra x^c$ denote the complex
conjugation on $\g$ with respect to the real form $\g_0= Lie(G)$.
A root $\alpha$ is called {\it imaginary} if
$\alpha^c=-\alpha$ and it is called {\it real} if
$\alpha^c=\alpha$. Every root space $\g_\alpha$ is one dimensional
and has a generator $X_\alpha$ satisfying:
$$
[X_\alpha,X_{-\alpha}]=Y_\alpha\ \ \ \ {\rm with}\ \ \ \
\alpha(.)=B(Y_\alpha,.)
$$ $$
B(X_\alpha,X_{-\alpha})=1
$$
and $X_\alpha^c=X_{\alpha^c}$ if $\alpha$ is non-imaginary and
$X_\alpha^c=\pm X_{-\alpha}$ if $\alpha$ is imaginary. An
imaginary root $\alpha$ is called {\it compact} if $X_\alpha^c=-X_{-\alpha}$ and {\it noncompact} otherwise.
Let $\Phi_n$ be the set of noncompact imaginary roots and choose a
set $\Phi^+$ of positive roots such that for $\alpha\in\Phi^+$
nonimaginary we have that $\alpha^c\in\Phi^+$. Let $W=W(G,H)$ be
the {\it Weyl group} of $(G,H)$, that
is
$$
W\= \frac{{\rm normalizer}(H)}{{\rm centralizer}(H)}.
$$
Let ${\rm rk}_\R(G)$ be the dimension of a maximal $\R$-split
torus in $G$ and let $\nu =\dim G/K -{\rm rk}_\R(G)$. We define
the {\it Harish-Chandra constant} of $G$ by
$$
c_G\= (-1)^{|\Phi_n^+|} (2\pi)^{|\Phi^+|}
2^{\nu/2}\frac{v(T)}{v(K)}|W|.
$$

\subsection{Normalization of Haar measures}
Although the results will not depend on normalizations we will
need to normalize Haar measures for the computations along the
way. First for any compact subgroup $C\subset G$ we normalize its
Haar measure so that it has total mass one, i.e., $\vol(C)=1$.
Next let $H\subset G$ be a reductive subgroup, and let $\theta_H$
be a Cartan involution on $H$ with fixed point set $K_H$. The
same way as for $G$ itself the form $B$ restricted to the Lie
algebra of $H$ induces a Riemannian metric on the manifold
$X_H=H/K_H$. Let $dx$ denote the volume element of that metric.
We get a Haar measure on $H$ by defining
$$
\int_Hf(h) dh \= \int_{X_H}\int_{K_H} f(xk) dk dx
$$
for any continuous function of compact support $f$ on $H$.

Let $P\subset G$ be a {\it parabolic subgroup} of $G$ (see \cite{wall-rg1} 2.2). Let $P=MAN$ be the
{\it Langlands decomposition} of
$P$. Then $M$ and $A$ are reductive, so there Haar measures can be
normalized as above. Since $G=PK=MANK$ there is a unique Haar
measure $dn$ on the unipotent radical $N$ such that for any
constant function $f$ of compact support on $G$ it holds:
$$
\int_G f(x)dx \= \int_M\int_A\int_N\int_K f(mank) dkdndadm.
$$
Note that these normalizations coincide for Levi subgroups with
the ones met by Harish-Chandra in (\cite{HC-HA1} sect. 7).

\subsection{Invariant distributions}
In this section we shall throughout assume that $G$ is a real
reductive group of inner type. A distribution $T$ on $G$, i.e., a
continuous linear functional $T: C_c^\infty(G)\ra \C$ is called
{\it invariant} if for any $f\in
C_c^\infty(G)$ and any $y\in G$ it holds: $T(f^y)=T(f)$, where
$f^y(x)=f(yxy^{-1})$. Examples are:
\begin{itemize}
\item orbital integrals: $f\mapsto \CO_g(f)=\int_{G_g\bs G}f(x^{-1}gx)dx$ and
\item traces: $f\mapsto \tr\pi(f)$ for $\pi\in\hat{G}$.
\end{itemize}
These two examples can each be expressed in terms of the other.
Firstly, Harish-Chandra proved that for any $\pi\in\hat{G}$ there
exists a conjugation invariant locally integrable function
$\Theta_\pi$ on $G$ such that for any $f\in C_c^\infty(G)$
$$
\tr\pi(f)\=\int_Gf(x)\Theta_\pi(x)  dx.
$$
Recall the {\it Weyl integration formula} which says that for any
integrable function $\ph$ on $G$ we have
$$
\int_G\ph(x)dx\=
\sum_{j=1}^r\rez{|W(G,H_j)|}\int_{H_j}\int_{G/H_j}\ph(xhx^{-1})|\det(1-h|\g
/\h_j)|dx dh,
$$
where $H_1,\dots,H_r$ is a maximal set of nonconjugate Cartan
subgroups in $G$ and for each Cartan subgroup $H$ we let
$W(G,H)$ denote its Weyl group, i.e., the quotient
of the normalizer of $H$ in $G$ by its centralizer.

An element $x$ of $G$ is called {\it regular} if
its centralizer is a Cartan subgroup. The set of regular elements
$G^{reg}$ is open and dense in $G$ with complement of measure zero. Therefore the integral above
can be taken over $G^{reg}$ only. Letting $H_j^{reg}:=H_j\cap
G^{reg}$ we get

\begin{proposition}
Let $N$ be a natural number bigger than $\frac{\dim G}{2}$, then
for any $f\in L_{2N}^1(G)$ and any $\pi\in\hat{G}$ we have
$$
\tr\pi(f)\=\sum_{j=1}^r\rez{|W(G,H_j)|}\int_{H_j^{reg}}
\CO_h(f)\Theta_\pi(h)|\det(1-h|\g /\h_j)| dh.
$$
\end{proposition}

That is, we have expressed the trace distribution in terms of
orbital integrals. In the other direction it is also possible to
express semisimple orbital integrals in terms of traces.

At first let $H$ be a $\theta$-stable Cartan subgroup of $G$. Let
$\h$ be its complex Lie algebra and let $\Phi =\Phi(\g,\h)$ be the
set of roots. Let $x\ra x^c$ denote the complex conjugation on
$\g$ with respect to the real form $\g_0= Lie_\R(G)$. Choose an
ordering $\Phi^+\subset \Phi$ and let $\Phi^+_I$ be the set of
positive imaginary roots. To any root $\alpha\in\Phi$ let
\begin{eqnarray*}
H &\ra & \C^\times\\
 h &\mapsto & h^\alpha
\end{eqnarray*}
be its character, that is, for $X\in\g_\alpha$ the root space to
$\alpha$ and any $h\in H$ we have $\Ad(h)X=h^\alpha X$. Now put
$$
\ '\lap_I(h)\= \prod_{\alpha\in\Phi_I^+}(1-h^{-\alpha}).
$$
Let $H=AT$ where $A$ is the connected split component and $T$ is
compact. An element $at\in AT=H$ is called {\it split
regular} if the centralizer of $a$ in $G$
equals the centralizer of $A$ in $G$. The split regular elements
form a dense open subset containing the regular elements of $H$.
Choose a parabolic $P$ with split component $A$, so $P$ has
Langlands decomposition $P=MAN$. For $at\in AT =H$ let
\begin{eqnarray*}
\lap_+(at) &=& \left| \det((1-\Ad((at)^{-1}))|_{\g /\a\oplus
\m})\right|^{\rez{2}}\\ {}\\
 & =& \left|\det((1-\Ad((at)^{-1}))|_\n)\right| a^{\rho_P}\\ {}\\
 &=&
 \left|\prod_{\alpha\in\Phi^+-\Phi_I^+}(1-(at)^{-\alpha})\right|
 a^{\rho_P},
\end{eqnarray*}
where $\rho_P$ is the half of the sum of the roots
in $\Phi(P,A)$, i.e., $a^{2\rho_P}=\det(a|\n)$. We will also write
$h^{\rho_P}$ instead of $a^{\rho_P}$.

For any $h\in H^{reg} = H\cap G^{reg}$ let
$$
'F_f^H(h)\= 'F_f(h)\= \ '\lap_I(h) \lap_+(h) \int_{G/A}
f(xhx^{-1}) dx.
$$
It then follows directly from the definitions that for $h\in
H^{reg}$ it holds
$$
\CO_h(f) \= \frac{'F_f(h)}
                 {h^{\rho_P}\det(1-h^{-1}|(\g/\h)^+)},
$$
where $(\g/\h)^+$ is the sum of the root spaces attached to
positive roots. There is an extension of this identity to
nonregular elements as follows: For $h\in H$ let $G_h$ denote its
centralizer in $G$. Let $\Phi^+(\g_h,\h)$ be the set of positive
roots of $(\g_h,\h)$. Let
$$
\varpi_h \= \prod_{\alpha\in\Phi^+(\g_h,\h)}Y_\alpha,
$$
then $\varpi_h$ defines a left invariant differential operator on
$G$.

\begin{lemma}
For any $f\in L_{2N}^1(G)$ and any $h\in H$ we have
$$
\CO_h(f) \= \frac{\varpi_h 'F_f(h)}
                 {c_h h^{\rho_P}\det(1-h^{-1}|(\g/\g_h)^+)}.
$$
\end{lemma}

\prf This is proven in section 17 of \cite{HC-DS}. \qed

Our aim is to express orbital integrals in terms of traces of
representations. By the above lemma it is enough to express
$'F_f(h)$ it terms of traces of $f$ when $h\in H^{reg}$. For this
let $H_1=A_1T_1$ be another $\theta$-stable Cartan subgroup of $G$
and let $P_1=M_1A_1N_1$ be a parabolic with split component $A_1$.
Let $K_1=K\cap M_1$. Since $G$ is connected the compact group
$T_1$ is an abelian torus and its unitary dual $\widehat{T_1}$ is
a lattice. The Weyl group $W=W(M_1,T_1)$ acts on $\widehat{T_1}$
and $\widehat{t_1}\in\widehat{T_1}$ is called {\it regular} if its
stabilizer $W(\widehat{t_1})$ in $W$ is trivial. The regular set
$\widehat{T_1}^{reg}$ modulo the action of $W(K_1,T_1)\subset
W(M_1,T_1)$ parameterizes the discrete series representations of
$M_1$ (see \cite{Knapp}). For $\widehat{t_1}\in\widehat{T_1}$
Harish-Chandra \cite{HC-S} defined a distribution
$\Theta_{\widehat{t_1}}$ on $G$ which happens to be the trace of
the discrete series representation $\pi_{\widehat{t_1}}$ attached
to $\widehat{t_1}$ when $\widehat{t_1}$ is regular. When
$\widehat{t_1}$ is not regular the distribution
$\Theta_{\widehat{t_1}}$ can be expressed as a linear combination
of traces as follows. Choose an ordering of the roots of
$(M_1,T_1)$ and let $\Omega$ be the product of all positive roots.
For any $w\in W$ we have $w\Omega = \epsilon(w)\Omega$ for a
homomorphism $\epsilon : W\ra \{ \pm 1\}$. For nonregular
$\widehat{t_1}\in\widehat{T_1}$ we get
$\Theta_{\widehat{t_1}}=\rez{|W(\widehat{t_1})|}\sum_{w\in
W(\widehat{t_1})}\epsilon(w)\Theta'_{w,\widehat{t_1}}$, where
$\Theta'_{w,\widehat{t_1}}$ is the character of an irreducible
representation $\pi_{w,\widehat{t_1}}$ called a limit of discrete
series representation. We will write $\pi_{\widehat{t_1}}$ for the
virtual representation $\rez{|W(\widehat{t_1})|}\sum_{w\in
W_{\widehat{t_1}}}\epsilon(w)\pi_{w,\widehat{t_1}}$.

Let $\nu :a\mapsto a^\nu$ be a unitary character of $A_1$ then
$\widehat{h_1}=(\nu,\widehat{t_1})$ is a character of
$H_1=A_1T_1$. Let $\Theta_{\widehat{h_1}}$ be the character of the
representation $\pi_{\widehat{h_1}}$ induced parabolically from
$(\nu,\pi_{\widehat{t_1}})$. Harish-Chandra has proven

\begin{theorem}\label{inv-orb-int}
Let $H_1,\dots,H_r$ be maximal a set of nonconjugate
$\theta$-stable Cartan subgroups with split components $A_1,\dots,A_r$. Let $H=H_j$ for some $j$ with split component $A$. Then
for each $j$ there exists a continuous function $\Phi_{H|H_j}$ on
$H^{reg}\times \hat{H_j}$ such that for $h\in H^{reg}$ it holds
$$
'F_f^H(h)\= \sum_{j=1}^r
\int_{\hat{H_j}}\Phi_{H|H_j}(h,\widehat{h_j})\
\tr\pi_{\widehat{h_j}}(f)\ d\widehat{h_j}.
$$
Further $\Phi_{H|H_j}=0$ unless there is $g\in G$ such that
$gAg^{-1}\subset A_j$. Finally for $H_j=H$ the function can be
given explicitly as
\begin{eqnarray*}
\Phi_{H|H}(h,\hat{h}) &=& \rez{|W(G,H)|}\sum_{w\in
W(G,H)}\epsilon(w|T)\langle w\hat{h},h\rangle\\
 &=& \rez{|W(G,H)|} \ '\lap(h)\Theta_{\hat{h}}(h),
\end{eqnarray*}
where $\ '\lap = \lap_+\ '\lap_I$.
\end{theorem}

\prf \cite{HC-S}. \qed

\subsection{Smoothness of induced functions}
Let $M$ be a smooth (i.e. $C^\infty$) manifold and let $D\subset M$ be an open subset. Let $S=M\setminus D$ be its complement.
A real or complex valued function $f$ on $D$ is said to \emph{vanish to order at least $k\in\Z$} at a point $s\in S$, if there exists an open neighbourhood $U$ of $s$ in $M$ such that the function
$$
u\ \mapsto\ \frac{|f(u)|}{d(u,s)^k}
$$
is bounded above on $U\cap D$.
Here $d(x,y)$ is the distance function attached to a Riemannian metric on $M$.
By its local nature, this notion does not depend on the choice of the metric.
Likewise, we say that $f$ \emph{vanishes to order at most $k$} at s if the function $u\mapsto\frac{|f(u)|}{d(u,s)^k}$ is bounded away from zero on $U\cap D$.
If both conditions hold with the same $k\in\Z$, we say that $f$ vanishes to order $k$ or that $f$ has order $k$ at $s$.
If $k$ is negative, we then also say that $f$ \emph{has a pole of order $-k$} at $s$.
Finally, we say that the function $f$ vanishes to order at most/at least $k$ on $S$ if it does so for every $s\in S$.

Let $M'$ be another smooth manifold and let $F\colon M\to M'$ be a smooth map.
Let $E\to M$ and $E'\to M'$ be smooth vector bundles and $\tilde F\colon E\to E'$ be a smooth linear lift of $F$, i.e., $\tilde F$ is smooth and maps the vector space $E_m$ linearly to the vector space $E'_{F(m)}$ for every $m\in M$.
Then we say that $\tilde F$ \emph{vanishes to order at least/at most $k$ at $S$} if for every two sections $\sigma$ of $E$ and $\al$ of $(E')^*$ the function
$$
m\ \mapsto\ \al(\tilde F(\sigma(m)))
$$
does so.
We speak of any such function as an \emph{entry} of $\tilde F$.
An example of a smooth linear lift is the differential,
$$
TF\colon TM\ \to\ TM'.
$$

An open subset $C$ of $D$ is called \emph{full}, if the boundary of $C$ is contained in $S$.
This is equivalent to saying that $C$ is a union of connected components of $D$.

For $x\in \R$ let $[x]$ be the largest integer with $[x]\le x$.

\begin{proposition}\label{smoothness}
Let $C$ be a full subset of $D$.
Let $F\colon M\to M'$ be smooth and assume that $F$ restricted to ${C}$ is a diffeomorphism with open image and that the bounday of $F(C)$ is contained in $F(S)$.
Assume further that $\det(TF)$ vanishes to order at most $k\in\N$ at $S$.
Let $f$ be a real or complex valued function on $C$ that vanishes to order at least $j\in\N$ at $S$ and is $j$-times continuously differentiable inside $C$.
Extend $f$ to a continuous function on $M$ by setting $f\equiv 0$ outside $C$.
Assume that $f$ factors over $F$, so $f$ induces a function $f'$ on $M'$ which vanishes outside $F(C)$.

Then the function $f'$ is at least $r$-times continuously differentiable on $M'$, where
$$
r\= \left[\frac{j-1}{r+1}\right].
$$
\end{proposition}

\prf
Let $m\in M$.
Taking local co-ordinates on $M$ and $M'$ the tangential $TF$ can be viewed as a Jacobi matrix $JF$.
The chain rule implies that on the open set $F(C)$ one has
$$
T(F^{-1})\= (TF)^{-1}\= \frac 1{\det(JF)}(JF)^\#,
$$
where for a matrix $A$ we write $A^\#$ for its complementary matrix.
Note that the entries of $A^\#$ are polynomials in the entries of $A$.
Since $\det(TF)$ vanishes to order at most $k$ at $S$, it follows that $T(F^{-1})$, which is defined on $F(C)$, has a pole of order at most $k$ at the boundary of $F(C)$.
On the open set $F(C)$ we can write $f'= f\circ F^{-1}$, so $Tf'=Tf\circ TF^{-1}$.
Since $f$ vanishes at $S$ to order at least $j$ it follows that $Tf$ vanishes at $S$ to order at least $j-1$, so $Tf'$ vanishes at $S'$ to order at least $j-1-k$.

We now finish the proof of Proposition \ref{smoothness} by induction of $j$.
First assume $j\le r+1$.
Then, since $j\ge 1$, the function $f'$ extends to a continuous function on $M'$ and so $f'\in C^0(M')$ as claimed.

For the induction step assume $j>r+1$.
Then $Tf'$ vanishes on $S'$ to order at least $j-1-k\in \N$.
This means that every entry of $Tf'$ vanishes at least to this order.
We pick an entry $g'$ and consider the function $g=g'\circ F$ on $M$.
Then $g$ vanishes at least to order $j-1-k$ at $S$ and by induction hypothesis the function $g'$ is continuously differentiable up to order $\left[ \frac{j-1-k-1}{k+1}\right]=\left[\frac{j-1}{k+1}\right] -1$.
Since $g'$ is an arbitrary entry of $Tf'$, it follows that $f'$ is of class $C^r$ with $r=\left[\frac{j-1}{k+1}\right]$.
\qed

\section{Euler-Poincar\'e functions}

\label{ep} In this section we generalize the construction
of pseudo-coefficients and Euler-Poincar\'e
functions \cite{CloDel, Lab} to non-connected groups. Here $G$ will be a real
reductive group. It will be assumed that $G$ has a
compact Cartan subgroup. It then follows that $G$ has
compact center.

\subsection{Existence}
Fix a maximal compact subgroup $K$ of $G$ and a Cartan $T$ of $G$
which lies inside  $K$. The group $G$ is called {\it orientation
preserving} if $G$ acts by
orientation preserving diffeomorphisms on the manifold $X=G/K$.
For example, the group $G=SL_2(\R)$ is orientation preserving but
the group $PGL_2(\R)$ is not. Recall the Cartan decomposition
$\g_0=\k_0\oplus\p_0$. Note that $G$ is orientation preserving if
and only if its maximal compact subgroup $K$ preserves
orientations on $\p_0$.

\begin{lemma} \label{orient}
The following holds:
\begin{itemize}
\item
Any connected group is orientation preserving.
\item
If $X$ carries the structure of a complex manifold which is left
stable by $G$, then $G$ is orientation preserving.
\end{itemize}
\end{lemma}

\prf The first is clear. The second follows from the fact that
biholomorphic maps are orientation preserving.
 \qed

Let $\t$ be the complexified Lie algebra of the Cartan
subgroup $T$. We choose an ordering of the roots $\Phi(\g
,\t)$ of the pair $(\g ,\t)$
\cite{wall-rg1}. This choice induces a decomposition $\p = \p_-
\oplus \p_+$, where $\p_\pm$ is the sum of the
positive/negative root spaces which lie in $\p$. As usual denote
by $\rho$ the half sum of the positive roots. The
chosen ordering induces an ordering of the {\it compact
roots} $\Phi(\k ,\t)$ which form a subset of
the set of all roots $\Phi(\g,\t)$. Let $\rho_K$ denote the half
sum of the positive compact roots. Recall that a function $f$ on
$G$ is called {\it $K$-central} if
$f(kxk^{-1})=f(x)$ for all $x\in G$, $k\in K$. For any
$K$-representation $(\rho,V)$ let $V^K$ denote the space of
$K$-fixed vectors, i.e., 
$$
V^K \=\{ v\in V | \rho(k)v=v\ \forall k\in K\}.
$$
Let $(\tau,V_\tau)$ be a representation of $K$ on a finite
dimensional complex vector space $V_\tau$. Let $V_\tau^*$ be the
dual space then there is a representation $\breve{\tau}$ on
$V_{\breve{\tau}}:=V_\tau^*$ given by
$$
\breve{\tau}(k)\alpha(v)\= \alpha(\tau(k^{-1})v),
$$
for $k\in K$, $\alpha\in V_\tau^*$ and $v\in V_\tau$. This
representation is called the {\it contragredient} or {\it dual}
representation. 
The restriction from $G$ to $K$ gives a ring
homomorphism of the representation rings:
$$
res_K^G:\Rep(G)\ra\Rep(K).
$$

\begin{theorem}\label{exist-ep}
(Euler-Poincar\'e functions) Let $(\tau ,V_\tau)$ a finite
dimensional representation of $K$. If $G$ is orientation
preserving or $\tau$ lies in the image of $res^G_K$, then there
is a compactly supported smooth $K$-central $K$-finite function $f_\tau$ on
$G$ such that for every admissible representation $(\pi
,V_\pi)$ of $G$ we have
$$
\tr\ \pi (f_\tau) \= \sum_{p=0}^{\dim (\p)} (-1)^p \dim (V_\pi
\otimes \wedge^p\p \otimes V_{\breve{\tau}})^K.
$$
We call $f_\tau$ an {\it Euler-Poincar\'e function} for $\tau$.
Note that, since $f$ is $K$-finite and $\pi$ is admissible, the operator $\pi(f)$ has finite rank, so the trace exists.

If $G$ is orientation preserving and $K$ leaves invariant the decomposition
$\p=\p_+\oplus\p_-$ then there is a compactly supported smooth
$K$-central function $g_\tau$ on $G$ such that for every
admissible representation $(\pi ,V_\pi)$ we have
$$
\tr\ \pi (g_\tau) \= \sum_{p=0}^{\dim (\p_-)} (-1)^p \dim (V_\pi
\otimes \wedge^p\p_- \otimes V_{\breve{\tau}})^K.
$$
\end{theorem}

If the representation $\tau$ lies in the image of
$res_K^G$ and the group $G$ is connected then the theorem is well
known, \cite{CloDel}, \cite{Lab}. 

\vspace{10pt}

\prf  We start with the case when $G$ is
orientation preserving. Without loss of generality assume $\tau$
irreducible. Suppose given a function $f$ which
satisfies the claims of the theorem except that it is not
necessarily $K$-central, then the function
$$
x\mapsto \int_K f(kxk^{-1})dk
$$
will satisfy all claims of the theorem. Thus one only needs to
construct a function having the claimed traces.

If $G$ is orientation preserving the adjoint action gives a
homomorphism $K\ra SO(\p)$. If this homomorphism happens to lift
to the double cover $Spin(\p)$ \cite{LawMich} we let $\tilde{G}=G$
and $\tilde{K}=K$. In the other case we apply the

\begin{lemma}\label{double-cover} If the homomorphism $K\ra
SO(\p)$ does not factor over the spin group $Spin(\p)$ then there is a double
covering $\tilde{G}\ra G$ such that with $\tilde{K}$ denoting the
inverse image of $K$ the induced homomorphism $\tilde{K}\ra
SO(\p)$ factors over $Spin(\p)\ra SO(\p)$. Moreover the kernel of
the map $\tilde{G}\ra G$ lies in the center of $\tilde{G}$
\end{lemma}

\prf At first $\tilde{K}$ is given by the pullback diagram:
$$
\begin{array}{ccc}
\tilde{K} & \ra & Spin(\p)\\
\downarrow & {} & \downarrow \alpha\\
K & \begin{array}{c}\Ad\\ \ra\\ {}\end{array} & SO(\p)
\end{array}
$$
that is, $\tilde{K}$ is given as the set of all $(k,g)\in K\times
Spin(\p)$ such that $\Ad(k)=\alpha(g)$. Then $\tilde{K}$ is a
double cover of $K$.

Next we use the fact that $K$ is a retract of $G$ to show that the
covering $\tilde{K}\ra K$ lifts to $G$. Explicitly let
$P=\exp(\p_0)$ then the map $K\times P\ra G, (k,p)\mapsto kp$ is a
diffeomorphism \cite{wall-rg1}. Let $g\mapsto
(\underline{k}(g),\underline{p}(g))$ be its inverse map. We let
$\tilde{G}\= \tilde{K}\times P$ then the covering $\tilde{K}\ra K$
defines a double covering $\beta :\tilde{G}\ra G$. We have to
install a group structure on $\tilde{G}$ which makes $\beta$ a
homomorphism and reduces to the known one on $\tilde{K}$. Now let
$k,k'\in K$ and $p,p'\in P$ then by
$$
k'p'kp\= k'k\ k^{-1}p'kp
$$
it follows that there are unique maps $a_K : P\times P\ra K$ and
$a_P :P\times P\ra P$ such that
\begin{eqnarray*}
\underline{k}(k'p'kp)&=& k'k a_K(k^{-1}p'k,p)\\
\underline{p}(k'p'kp)&=& a_P(k^{-1}p'k,p).
\end{eqnarray*}
Since $P$ is simply connected the map $a_K$ lifts to a map
$\tilde{a}_K: P\times P\ra \tilde{K}$. Since $P$ is connected
there is exactly one such lifting with $\tilde{a}_K(1,1)=1$. For $k\in\tilde K$ let $\bar k$ be its image in $K$. Now
the map
\begin{eqnarray*}
(\tilde{K}\times P)\times(\tilde{K}\times P)&\ra& \tilde{K}\times P\\
(k',p'),(k,p) &\mapsto&
(kk'\tilde{a}_K(\bar k^{-1}p'\bar k,p),a_P(\bar k^{-1}p'\bar k,p))
\end{eqnarray*}
defines a multiplication on $\tilde{G}=\tilde{K}\times P$ with
the desired properties.

Finally $\ker(\beta)$ will automatically be central because it is
a normal subgroup of order two.
\qed

Let $S$ be the spin representation of $Spin(\p)$ (see
\cite{LawMich}, p.36). It splits as a direct sum of two distinct
irreducible representations
$$
S\= S^+\oplus S^-.
$$
We will make use of the following properties of the spin representation.
\begin{itemize}
\item
The virtual representation
$$
(S^+- S^-)\otimes (S^+- S^-)
$$
is isomorphic to the adjoint representation on $\wedge^{even}\p
-\wedge^{odd}\p$ (see \cite{LawMich}, p. 36).
\item
If $K$ leaves invariant the spaces $\p_-$ and $\p_+$, as is the
case when $X$ carries a holomorphic structure fixed by $G$, then
there is a one dimensional representation $\epsilon$ of
$\tilde{K}$ such that
$$
(S^+-S^-)\otimes\epsilon \ \cong\
\wedge^{even}\p_--\wedge^{odd}\p_-.
$$
\end{itemize}
The proof of this latter property will be given in section
\ref{appA}.

\begin{theorem} \label{existh}
(Pseudo-coefficients) Assume that the group $G$ is
orientation preserving. Then for any finite dimensional
representation $(\tau ,V_\tau)$ of $\tilde{K}$ there is a
compactly supported smooth function $h_\tau$ on
$\tilde{G}$ such that for every admissible
representation $(\pi ,V_\pi)$ of $\tilde{G}$,
$$
\tr\ \pi (h_\tau) = \dim (V_\pi \otimes S^+ \otimes
V_{\breve{\tau}})^{\tilde{K}} - \dim (V_\pi \otimes S^- \otimes
V_{\breve{\tau}})^{\tilde{K}}.
$$
\end{theorem}

The functions given in this theorem are also known as
pseudo-coefficients \cite{Lab}. This result generalizes
the one in \cite{Lab} in several ways. First, the group
$G$ needn't be connected and secondly the representation
$\tau$ needn't be spinorial. The proof of this theorem
relies on the following lemma.

\begin{lemma}\label{2.5}
Let $(\pi ,V_\pi)$ be an irreducible admissble representation of
$\tilde{G}$ and assume
$$
\dim (V_\pi \otimes S^+ \otimes V_{\breve{\tau}})^{\tilde{K}} -
\dim (V_\pi \otimes S^- \otimes V_{\breve{\tau}})^{\tilde{K}}
\neq 0,
$$
then the Casimir eigenvalue satisfies $\pi (C) = \breve{\tau}
(C_K) - B(\rho)+B(\rho_K)$.
\end{lemma}

\prf Let the $\tilde{K}$-invariant operator
 $$
 d_\pm : V_\pi\otimes S^\pm \ra V_\pi\otimes S^\mp
 $$
be defined by
 $$
 d_\pm : v\otimes s \mapsto \sum_i\pi(X_i)v\otimes c(X_i)s,
 $$
where $(X_i)$ is an orthonormal base of $\p$. The formula of
Parthasarathy, \cite{AtSch}, p. 55 now says
 $$
 d_-d_+\= d_+d_- \= \pi\otimes s^\pm(C_K) -\pi(C)\otimes 1
 -B(\rho)+B(\rho_K),
 $$
 where $s^\pm$ is the representation on $S^\pm$.
Our assumption leads to $ker(d_+ d_-) \cap \pi \otimes S(\tau)
\neq 0$, where $\pi\otimes S(\tau)$ is the $\tau$ K-type of $\pi\otimes S$. 
But on this space the $K$-Casimir $C_K$ acts by the scalar $\tau(C_K)$, so that we get $0=\tau(C_K)-\pi(C) -B(\rho) +B(\rho_K)$.
\qed

For the proof of Theorem \ref{existh} let $(\tau,V_\tau)$ a finite
dimensional irreducible unitary representation of $\tilde{K}$ and
write $E_\tau$ for the $\tilde{G}$-homogeneous vector bundle over
${X}=\tilde{G}/\tilde{K}$ defined by $\tau$. The space of smooth
sections $\Ga^\infty (E_\tau)$ may be written as $\Ga^\infty
(E_\tau) = (C^\infty (\tilde{G}) \otimes V_\tau)^{\tilde{K}}$,
where $\tilde{K}$ acts on $C^\infty (\tilde{G})$ by right
translations. The Casimir operator $C$ of $\tilde{G}$ acts on
this space and defines a second order differential operator
$C_\tau$ on $E_\tau$. On the space of $L^2$-sections
$L^2(X,E_\tau)=(L^2(\tilde{G}) \otimes V_\tau)^{\tilde{K}}$ this
operator is formally selfadjoint with domain, say, the compactly
supported smooth functions and extends to a selfadjoint operator.
Let $g$ be a Paley-Wiener function on $\R$, i.e., $g$ is the Fourier-transform of a smooth function of compact support.
Then $g$ extends to a holomorphic function on $\C$.
Assume that $g$ is even, i.e., $g(z)=g(-z)$.
Then the power series of $g(z)$ around zero contains only even powers of $z$.
So there is an entire function $f$ such that $g(z)=f(z^2)$.
Then $f|_\R$ is a Schwartz-Bruhat function.
In \cite{mult} Lemma 2.9 and Lemma 2.11 it is shown that there exists a smooth function of compact support $\tilde f_\tau$ on $G$ such that for every irreducible unitary representation 
 $\pi$ of ${\tilde{G}}$:
\begin{equation}\label{ftau}
\pi(\tilde{f}_\tau) = f(\pi(C))\Pr_{\breve\tau},
\end{equation}
where $\Pr_{\breve\tau}$ is the orthogonal projection to the $K$-type $\breve\tau$, and $\pi(C)$ denotes the Casimir
eigenvalue on $\pi$. 
This formula holds as well for $\pi$ being an irreducible admissible representation, as is seen as follows.
First let $\pi=\pi_{\xi,\la}$ be a representation induced from a minimal parabolic $P=MAN$, where $\xi\in\hat M$ and $\la\in\a^*$.
If $\la$ is imaginary and generic, then $\pi$ is irreducible unitary, so one has
$$
\pi_{\xi,\la}(\tilde{f}_\tau) = f(\pi_{\xi,\la}(C))\Pr_{\breve\tau}.
$$ 
Both sides of this equality are holomorphic functions in $\la$ with values in the finite dimensional space $\End(\pi(\breve\tau))$, hence the result holds for any $\la$ by the identity theorem for holomorphic functions.
Finally, any irreducible admissible $\pi$ is a sub-representation of an induced representation, hence the formula (\ref{ftau}) holds for every irreducible admissible representation $\pi$.

For $\tau=\tau_1\oplus\tau_2$ let $\tilde f_\tau=\tilde f_{\tau_1}+\tilde f_{\tau_2}$, so $\tilde f_\tau$ is defined for all finite dimensional representations of $K$.
Next for a virtual representation $\tau=\tau_1-\tau_2$ we let $\tilde f_{\tau}=\tilde f_{\tau_1}-\tilde f_{\tau_2}$.

Choose $f$ as above such that $f(\breve{\tau} (C_K) -B(\rho)+B(\rho_K))=1$.
Such an $f$ clearly exists. 
Let $\tau\in\hat K$ and let $\ga$ be the virtual
representation of $\tilde{K}$ on the space
$$
V_\ga = (S^+-S^-)\otimes V_\tau,
$$
then set $h_\tau =\tilde{f}_\ga$. Then for an irreducible admissible representation $\pi$ of $G$ one has
$$
\tr\pi(h_\tau)\= f(\pi(C))\left(\dim (V_\pi\otimes S^+\otimes V_{\breve\tau})^{\tilde K}-\dim(V_\pi\otimes S^-\otimes V_{\breve\tau})^{\tilde K}\right).
$$
By Lemma \ref{2.5} and the choice of $f$ this gives
$$
\tr\pi(h_\tau)\= \left(\dim (V_\pi\otimes S^+\otimes V_{\breve\tau})^{\tilde K}-\dim(V_\pi\otimes S^-\otimes V_{\breve\tau})^{\tilde K}\right).
$$

So the function $h_\tau$ has the property claimed in Theorem \ref{existh} for irreducible $\pi$.
It immediately follows for direct sums of irreducibles.
Since the assertion only involves the trace of $\pi(h_\tau)$, it follows for arbitrary $\pi$ since it is valid for the semisimplification of $\pi$.
\qed

To get the first part of Theorem \ref{exist-ep} from Theorem
\ref{existh} one replaces $\tau$ in the proposition by the virtual
representation on $(S^+-S^-)\otimes V_\tau$. Since
$(S^+-S^-)\otimes (S^+-S^-)$ is as $\tilde{K}$ module isomorphic
to $\wedge^*\p$ we get the desired function, say $j$ on the group
$\tilde{G}$. Now if $\tilde{G}\ne G$ let $z$ be the nontrivial
element in the kernel of the isogeny $\tilde{G}\ra G$, then the
function
$$
f(x) \= \rez{2}(j(x)+j(zx))
$$
factors over $G$ and satisfies the claim.

To get the second part of the theorem one proceeds similarly
replacing $\tau$ by $\epsilon\otimes \tau$.

It remains to consider the case when $G$ is not orientation preserving, but $\tau$ lies in the image of the restriction map.
For this it suffices to show the claim in the case when $\tau$ is replaced by a finite dimensional irreducible representation $(\sigma,V_\sigma)$ of $G$.
Then one proceeds as in the proof of Theorem \ref{existh}, except that the role of Lemma \ref{2.5} is taken up by the following lemma.

\begin{lemma} 
Let $(\sigma ,V_\sigma)$ be an irreducible finite dimensional
representation of $G$. Let $(\pi ,V_\pi)$ be an irreducible
unitary representation of $G$ and assume
$$
\sum_{p=0}^{\dim (\p)} (-1)^p \dim (V_\pi \otimes \wedge^p\p
\otimes V_{\sigma})^{K} \neq 0,
$$
then the Casimir eigenvalues satisfy
 $$
 \pi(C) = \sigma(C).
 $$
\end{lemma}

\prf Recall that the Killing form of $G$ defines a
$K$-isomorphism between $\p$ and its dual $\p^*$, hence in
the assumption of the lemma we may replace $\p$ by $\p^*$. Let
$\pi_K$ denote the $(\g ,K)$-module of $K$-finite vectors in
$V_\pi$ and let $C^q(\pi_K\otimes V_\sigma) =
\Hom_{\k}(\wedge^q\p ,\pi_K\otimes V_\sigma) =
(\wedge^q\p^*\otimes\pi_K\otimes V_\sigma)^{\k}$ the standard
complex for the relative Lie algebra cohomology $H^q(\g ,\k
,\pi_K\otimes V_\sigma)$. Further
$(\wedge^q\p^*\otimes\pi_K\otimes V_\sigma)^{K_M}$ forms the
standard complex for the relative $(\g ,K)$-cohomology $H^q(\g,K ,\pi_K\otimes V_\sigma)$. In \cite{BorWall}, p.28 it is shown
that
$$
H^q(\g ,K ,\pi_K\otimes V_\sigma)=H^q(\g ,\k
,\pi_K\otimes V_\sigma)^{K/K^0}.
$$
Our assumption
implies $\sum_q(-1)^q \dim H^q(\g ,K ,\pi_K\otimes
V_\sigma) \ne 0$, therefore there is a $q$ with $0\ne
H^q(\g ,K ,\pi_K\otimes V_\sigma)=H^q(\g ,\k
,\pi_K\otimes V_\sigma)^{K/K^0}$, hence $H^q(\g ,\k
,\pi_K\otimes V_\sigma)\ne 0$. Now Proposition 3.1 on page
52 of \cite{BorWall} says that $\pi(C)\ne\sigma(C)$
implies that $H^q(\g,\k,\pi_K\otimes V_\sigma)=0$ for
all $q$. The claim follows. \qed

\subsection{Clifford algebras and Spin groups}\label{appA}
This section is solely given to provide a proof of the properties
of the spin representation used in the last section. We will
therefore not strive for the utmost generality but plainly state
things in the form needed. For more details the reader is
referred to \cite{LawMich}.

Let $V$ be a finite dimensional complex vector space and let
$q:V\ra\C$ be a non-degenerate quadratic form. We use the same
letter for the symmetric bilinear form:
$$
q(x,y) \= \rez{2}(q(x+y)-q(x)-q(y)).
$$
Let $SO(q)\subset GL(V)$ be the special orthogonal group of $q$.
The {\it Clifford algebra}  $Cl(q)$ will
be the quotient of the tensorial algebra
$$
TV \= \C\oplus V\oplus (V\otimes V)\oplus\dots
$$
by the two-sided ideal generated by all elements of the form
$v\otimes v+q(v)$, where $v\in V$.

This ideal is not homogeneous with respect to the natural
$\Z$-grading of $TV$, but it is homogeneous with respect to the
induced $\Z/2\Z$-grading given by the even and odd degrees. Hence
the latter is inherited by $Cl(q)$:
$$
Cl(q)\= Cl^0(q)\oplus Cl^1(q).
$$
For any $v\in V$ we have in $Cl(q)$ that $v^2=-q(v)$ and therefore
$v$ is invertible in $Cl(q)$ if $q(v)\ne 0$. Let $Cl(q)^\times$
be the group of invertible elements in $Cl(q)$. The algebra
$Cl(q)$ has the following universal property: For any linear map
$\ph :V\ra A$ to a $\C$-algebra $A$ such that $\ph(v)^2=-q(v)$
for all $v\in V$ there is a unique algebra homomorphism $Cl(v)\ra
A$ extending $\ph$.

Let $Pin(q)$ be the subgroup of the group $Cl(q)^\times$
generated by all elements $v$ of $V$ with $q(v)=\pm 1$. Let the
{\it complex spin group} be defined by
$$
Spin(q) \= Pin(q)\cap Cl^0(q),
$$
i.e., the subgroup of $Pin(q)$ of those elements which are
representable by an even number of factors of the form $v$ or
$v^{-1}$ with $v\in V$. Then $Spin(q)$ acts on $V$ by $x.v =
xvx^{-1}$ and this gives a fourfold covering: $Spin(q)\ra SO(q)$.

Assume the dimension of $V$ is even and let
$$
V \= V^+\oplus V^-
$$
be a {\it polarization},  that is
$q(V^+)=q(V^-)=0$. Over $\C$ polarizations always exist for even
dimensional spaces. By the nondegeneracy of $q$ it follows that
to any $v\in V^+$ there is a $\hat{v}\in V^-$ such that
$q(v,\hat{v})=-1$. Further, let $V^{-,v}$ be the space of all
$w\in V^-$ such that $q(v,w)=0$, then
$$
V^- \= \C \hat{v}\oplus V^{-,v}.
$$
Let
$$
S\= \wedge^* V^- \= \C\oplus V^-\oplus \wedge^2
V^-\oplus\dots\oplus\wedge^{top}V^-,
$$
then we define an action of $Cl(q)$ on $S$ in the following way:
\begin{itemize}
\item
for $v\in V^-$ and $s\in S$ let
$$
v.s\= v\wedge s,
$$
\item
for $v\in V^+$ and $s\in \wedge^*V^{-,v}$ let
$$
v.s\= 0,
$$
\item
and for $v\in V^+$ and $s\in S$ of the form $s=\hat{v}\wedge s'$
with $s'\in \wedge^*V^{-,v}$ let
$$
v.s\= s'.
$$
\end{itemize}
By the universal property of $Cl(V)$ this extends to an action of
$Cl(q)$. The module $S$ is called the {\it spin
module}. The induced action of $Spin(q)$
leaves invariant the subspaces
$$
S^+\=\wedge^{even}V^-,\ \ \ S^-\=\wedge^{odd}V^-,
$$
the representation of $Spin(q)$ on these spaces are called the
{\it half spin representations}.
Let $SO(q)^+$ the subgroup of all elements in $SO(q)$ that leave
stable the decomposition $V=V^+\oplus V^-$. This is a connected
reductive group isomorphic to $GL(V^+)$, since, let $g\in
GL(V^+)$ and define $\hat{g}\in GL(V^-)$ to be the inverse of the
transpose of $g$ by the pairing induced by $q$ then the map
$Gl(v)\ra SO(q)^+$ given by $g\mapsto (g,\hat{g})$ is an
isomorphism. In other words, choosing a basis on $V^+$ and a the
dual basis on $V^-$ we get that $q$ is given in that basis by
$\matrix{0}{\1}{\1}{0}$. Then $SO(q)^+$ is the image of the
embedding
\begin{eqnarray*}
GL(V^-) &\hookrightarrow& SO(q)\\
A &\mapsto& \matrix{A}{0}{0}{^tA^{-1}}.
\end{eqnarray*}

Let $Spin(q)^+$ be the inverse image of $SO(q)^+$ in $Spin(q)$.
Then the covering $Spin(q)^+\ra SO(q)^+\cong GL(V^-)$ is the
``square root of the determinant'', i.e., it is isomorphic to the
covering $\tilde{GL}(V^-)\ra GL(V^-)$ given by the pullback
diagram of linear algebraic groups:
$$
\begin{array}{ccc}
\tilde{GL}(V^-)     & \ra   & GL(1)\\
\downarrow         & {} &\downarrow x\mapsto x^2\\
GL(V^-)         & \begin{array}{c}\det\\ \ra\\ {}\end{array}&
GL(1).
\end{array}
$$
As a set, $\tilde{GL}(V^-)$ is given  as the set of all pairs
$(g,z)\in GL(V^-)\times GL(1)$ such that $\det(g)=z^2$ and the
maps to $GL(V^-)$ and $GL(1)$ are the respective projections.

\begin{lemma}\label{epsilon}
There is a one dimensional representation $\epsilon$ of
$Spin(q)^+$ such that
$$
S^\pm\otimes\epsilon\ \cong\ \wedge^\pm V^+
$$
as $Spin(q)^+$-modules, where $\wedge^\pm$ means the even or odd
powers respectively.
\end{lemma}

\prf Since $Spin(q)^+$ is a connected reductive group over $\C$
we can apply highest weight theory. If the weights of the
representation of $Spin(q)^+$ on $V$ are given by
$\pm\mu_1,\dots,\pm\mu_m$, then the weights of the half spin
representations are given by
$$
\rez{2}(\pm\mu_1\pm\dots\pm\mu_m)
$$
with an even number of minus signs in the one and an odd number
in the other case. Let $\epsilon = \rez{2}(\mu_1+\dots +\mu_m)$
then $\epsilon$ is a weight for $Spin(q)^+$ and $2\epsilon$ is
the weight of, say, the one dimensional representation on
$\wedge^{top}V^+$. By Weyl's dimension formula this means that
$2\epsilon$ is invariant under the Weyl group and therefore
$\epsilon$ is. Again by Weyl's dimension formula it follows that
the representation with highest weight $\epsilon$ is one
dimensional. Now it follows that $S^+\otimes\epsilon$ has the
same weights as the representation on $\wedge^+ V^+$, hence must be
isomorphic to the latter. The case of the minus sign is analogous.
\qed

\subsection{Orbital integrals}
It now will be shown that $\tr\pi(f_\tau)$ vanishes for a
principal series representation $\pi$. To this end let $P=MAN$ be
a nontrivial parabolic subgroup with $A\subset \exp(\p_0)$. Let
$(\xi ,V_\xi)$ denote an irreducible unitary representation of $M$
and $e^\nu$ a quasicharacter of $A$. Let $\pi_{\xi ,\nu}:= {\rm
Ind}_P^G (\xi \otimes e^{\nu}\otimes 1)$.

\begin{lemma} \label{pivonggleichnull}
We have $\tr\pi_{\xi ,\nu}(f_\tau) =0$.
\end{lemma}

\prf By Frobenius reciprocity we have for any irreducible unitary
representation $\ga$ of $K$:
$$
\Hom_K(\ga ,\pi_{\xi ,\nu}|_K) \cong \Hom_{K_M}(\ga |_{K_M},\xi ),
$$
where $K_M := K\cap M$. This implies that $\tr\pi_{\xi
,\nu}(f_\tau)$ does not depend on $\nu$. On the other hand
$\tr\pi_{\xi ,\nu}(f_\tau)\ne 0$ for some $\nu$ would imply
$\pi_{\xi ,\nu}(C) =\breve{\tau} (C_K) -B(\rho)+B(\rho_K)$ which
only can hold for $\nu$ in a set of measure zero.
\qed

Recall that an element $g$ of $G$ is called {\it
elliptic} if it lies in a compact Cartan subgroup.
Since the following relies on results of Harish-Chandra which were
proven under the assumption that $G$ is of inner type, we will
from now on assume this.

\begin{theorem} \label{orbitalint}
Assume that $G$ is of inner type. Let $g$ be a semisimple element
of the group $G$. If $g$ is not elliptic, then the orbital
integral $\O_g(f_\tau)$ vanishes. If $g$ is elliptic we may assume
$g\in T$, where $T$ is a Cartan in $K$ and then we have
$$
\O_g(f_\tau) \= \tr\ \tau(g)\ c_g^{-1}|W(\t ,\g_g)| \prod_{\alpha
\in \Phi_g^+}(\rho_g ,\alpha),
$$
where $c_g$ is Harish-Chandra's constant, it does only depend on
the centralizer $G_g$ of g. Its value is given in \ref{notations}.

\end{theorem}

\prf The vanishing of $\O_g(f_\tau)$ for nonelliptic semisimple
$g$ is immediate by the lemma above and Theorem \ref{inv-orb-int}.
So consider $g\in T\cap G'$, where $G'$ denotes the set of regular
elements. Note that for regular $g$ the claim is $\CO_g(f_\tau)
=\tr \tau(g)$. Assume the claim proven for regular elements, then
the general result follows by standard considerations as in
\cite{HC-DS}, p.32 ff. where however different Haar-measure
normalizations are used that produce a factor $[G_g:G_g^0]$,
therefore these standard considerations are now explained. Fix
$g\in T$ not necessarily regular. Let $y\in T^0$ be such that $gy$
is regular. Then
\begin{eqnarray*}
\tr \tau (gy) &=& \int_{T\bs G} f_\tau(x^{-1}gyx) dx\\
    &=&  \int_{T^0\bs G} f_\tau(x^{-1}gyx) dx\\
    &=&  \int_{G_g\bs G}\int_{T^0\bs G_g} f_\tau(x^{-1}z^{-1}gyzx)\ dz\ dx\\
    &=&  \int_{G_g\bs G}\sum_{\eta :G_g /G_g^0}
\rez{[G_g:G_g^0]}\int_{T^0\bs G_g^0}
f_\tau(x^{-1}\eta^{-1}z^{-1}gyz\eta x)\ dz\ dx.
\end{eqnarray*}
The factor $\rez{[G_g:G_g^0]}$ comes in by the Haar-measure
normalizations. On $G_g^0$ consider the function
$$
h(y) = f(x^{-1}\eta^{-1} yg\eta x).
$$
Now apply Harish-Chandra's operator $\omega_{G_g}$ to $h$ then
for the connected group $G_g^0$ it holds
$$
h(1) = \lim_{y\ra 1} c_{G_g^0}^{-1}\ \omega_{G_g^0}\ F_h^{G_g^0}(y),
$$
where $F_h$ is Harish-Chandra's invariant integral \cite{HC-DS}.
When $y$ tends to $1$ the $\eta$-conjugation drops out and the
claim follows.

So in order to prove the proposition one only has to consider the
regular orbital integrals. Next the proof will be reduced to the
case when the compact Cartan $T$ meets all connected components
of $G$. For this let $G^+=TG^0$ and assume the claim proven for
$G^+$. Let $x\in G$ then $xTx^{-1}$ again is a compact Cartan
subgroup. Since $G^0$ acts transitively on all compact Cartan
subalgebras it follows that $G^0$ acts transitively on the set of
all compact Cartan subgroups of $G$. It follows that there is a
$y\in G^0$ such that $xTx^{-1} =yTy^{-1}\subset TG^0 =G^+$, which
implies that $G^+$ is normal in $G$.

Let $\tau^+=\tau |_{G^+\cap K}$ and $f^+_{\tau^+}$ the
corresponding Euler-Poincar\'e function on $G^+$.

\begin{lemma}
$f^+_{\tau^+} = f_\tau |_{G^+}$
\end{lemma}

Since the Euler-Poincar\'e function is not uniquely determined the
claim reads that the right hand side is a EP-function for $G^+$.

\prf Let $\tau^+ =\tau |_{K^+}$, where $K^+ = TK^0 = K\cap G^+$.
Let $\ph^+\in(C_c^\infty(G^+)\otimes V_\tau)^{K^+}$, which may be
viewed as a function $\ph^+ : G^+ \ra V_\tau$ with $\ph^+(xk)
=\tau(k^{-1})\ph^+(x)$ for $x\in G^+$, $k\in K^+$. Extend $\ph^+$
to $\ph : G\ra V_\tau$ by $\ph(xk)=\tau(k^{-1})\ph^+(x)$ for
$x\in G^+$, $k\in K$. This defines an element of
$(C_c^\infty(G)\otimes V_\tau)^K$ with $\ph |_{G^+}=\ph^+$. Since
$C_\tau$ is a differential operator it follows $f(C_\tau)\ph
|_{G^+} =f(C_\tau)\ph^+$, so
$$
(\ph * \tilde{f}_\tau) |_{G^+} = \ph^+ * \tilde{f}_{\tau^+}.
$$
Considering the normalizations of Haar measures gives the lemma.
\qed

For $g\in T'$ we compute
\begin{eqnarray*}
\CO_g(f_\tau) &=& \int_{T\bs G} f_\tau (x^{-1} gx) dx\\
    &=& \sum_{y:G/G^+} \rez{[G:G^+]} \int_{y^{-1}Ty\bs G^+} f_\tau (x^{-1} y^{-1} g yx) dx,
\end{eqnarray*}
 where the factor $\rez{[G:G^+]}$ stems from normalization of
Haar measures and we have used the fact that $G^+$ is normal. The
latter equals
$$
\rez{[G:G^+]}\sum_{y:G/G^+}\CO^{G^+}_{y^{-1}gy}(f_\tau) =
\rez{[G:G^+]}\sum_{y:G/G^+}\CO^{G^+}_{y^{-1}gy}(f_{\tau^+}^+).
$$
Assuming the theorem proven for $G^+$, this is
$$
\rez{[G:G^+]}\sum_{y:G/G^+} \tr \tau(y^{-1}gy) =\tr\tau(g).
$$

From now on one thus may assume that the compact Cartan $T$ meets
all connected components of $G$. Let $(\pi,V_\pi)\in\hat{G}$.
Harish-Chandra has shown that for any $\ph\in C^\infty_c(G)$ the
operator $\pi(\ph)$ is of trace class and there is a locally
integrable conjugation invariant function $\Theta_\pi$ on $G$,
smooth on the regular set such that
$$
\tr \pi(\ph) = \int_G\ph(x)\Theta_\pi(x) dx.
$$
For any $\psi\in C^\infty(K)$ let
$\pi|_K(\psi)=\int_K\psi(k)\pi(k) dk$.

\begin{lemma} \label{character}
Assume $T$ meets all components of $G$. For any $\psi\in
C^\infty(K)$ the operator $\pi|_K(\psi)$ is of trace class and
for $\psi$ supported in the regular set $K' =K\cap G'$ we have
$$
\tr \pi|_K(\psi) = \int_K \psi(k) \Theta_\pi(k) dk.
$$
(For $G$ connected this assertion is in \cite{AtSch} p.16.)
\end{lemma}

\prf Let $V_\pi =\bigoplus_i V_\pi(i)$ be the decomposition of
$V_\pi$ into $K$-types. This is stable under $\pi|_K(\psi)$.
Harish-Chandra has proven $[\pi|_K :\tau]\le\dim \tau$ for any
$\tau\in\hat{K}$. Let $\psi =\sum_j\psi_j$ be the decomposition
of $\psi$ into $K$-bitypes. Since $\psi$ is smooth the sequence
$\norm{\psi_j}_1$ is rapidly decreasing for any enumeration of
the $K$-bitypes. Here $\norm{\psi}_1$ is the $L^1$-norm on $K$.
It follows that the sum $\sum_i\tr(\pi|_K(\psi)|V_\pi(i))$
converges absolutely, hence $\pi|_K(\psi)$ is of trace class.

Now let $S=\exp(\p_0)$ then $S$ is a smooth set of
representatives of $G/K$. Let $G.K=\cup_{g\in
G}gKg^{-1}=\cup_{s\in S}sKs^{-1}$, then, since $G$ has a compact
Cartan, the set $G.K$ has non-empty interior. Applying the Weyl
integration formula to $G$ and backwards to $K$ gives the
existence of a smooth measure $\mu$ on $S$ and a function $D$ with
$D(k)>0$ on the regular set such that
$$
\int_{G.K}\ph(x) dx = \int_S \int_K \ph(sks^{-1}) D(k) dk d\mu(s)
$$
for $\ph\in L^1(G.K)$. Now suppose $\ph\in C_c^\infty(G)$ with
support in the regular set. Then
\begin{eqnarray*}
\tr \pi(\ph) &=& \int_{G.K}\ph(x) \Theta_\pi(x) dx\\
    &=& \int_S \int_K \ph(sks^{-1}) D(k) \Theta_\pi(k)d\mu(s)\\
    &=& \int_K\int_S \ph^s(k) d\mu(s) D(k) \Theta_\pi(k) dk,
\end{eqnarray*}
where we have written $\ph^s(k)=\ph(sks^{-1})$. On the other hand
\begin{eqnarray*}
\tr\pi(\ph) &=& \tr \int_{G.K} \ph(x) \pi(x) dx\\
    &=& \tr \int_S\int_K \ph(sks^{-1})D(k) \pi(sks^{-1})dk d\mu(s)\\
    &=& \tr \int_S \pi(s) \pi|_K(\ph^s D)\pi(s)^{-1} d\mu(s)\\
    &=& \int_S \tr \pi |_K(\ph^s D) d\mu(s)\\
    &=& \tr\pi |_K\left( \int_S \ph^s d\mu(s) D\right).
\end{eqnarray*}
This implies the claim for all functions $\psi\in C_c^\infty(K)$
which are of the form
$$
\psi(k) = \int_S\ph(sks^{-1})d\mu(s) D(k)
$$
for some $\ph\in C_c^\infty(G)$ with support in  the regular set.
Consider the map
$$
\begin{array}{cccc}
F:& S\times K' & \ra & G.K'\\
{}& (s,k) & \mapsto & sks^{-1}
\end{array}
$$
Then the differential of $F$ is an isomorphism at any point and by
the inverse function theorem $F$ locally is a diffeomorphism. So
let $U\subset S$ and $W\subset K'$ be open sets such that
$F|_{U\times W}$ is a diffeomorphism. Then let $\alpha\in
C_c^\infty(U)$ and $\beta\in C_c^\infty(W)$, then define
$$
\phi(sks^{-1}) = \alpha(s)\beta(k)\ \ \ {\rm if}\ s\in U,\ k\in W
$$
and $\ph(g)=0$ if $g$ is not in $F(U\times W)$. We can choose the
function $\alpha$ such that $\int_S\alpha (s) d\mu(s) =1$. Then
$$
\int_S\ph(sks^{-1}) d\mu(s) D(k) = \beta(k) D(k).
$$
Since $\beta$ was arbitrary and $D(k) >0$ on $K'$ the lemma
follows.
\qed

Let $W$ denote the virtual $K$-representation on $\wedge^{even}\p
\otimes V_\tau - \wedge^{odd}\p\otimes V_\tau$ and write $\chi_W$
for its character.

\begin{lemma} \label{fin-lin-comb}
Assume $T$ meets all components of $G$, then for any
$\pi\in\hat{G}$ the function $\Theta_\pi\chi_W$ on $K'=K\cap G'$
equals a finite integer linear combination of $K$-characters.
\end{lemma}

\prf It suffices to show the assertion for $\tau =1$. Let $\ph$
be the homomorphism $K\ra O(\p)$ induced by the adjoint
representation, where the orthogonal group is formed with respect
to the Killing form. We claim that $\ph(K)\subset SO(\p)$, the
subgroup of elements of determinant one. Since we assume $K=K^0
T$ it suffices to show $\ph(T)\subset SO(\p)$. For this let $t\in
T$. Since $t$ centralizes $\t$ it fixes the  decomposition $\p
=\oplus_\alpha \p_\alpha$ into one dimensional root spaces. So
$t$ acts by a scalar, say $c$ on $\p_\alpha$ and by $d$ on
$\p_{-\alpha}$. There is $X\in\p_\alpha$ and $Y\in\p_{-\alpha}$
such that $B(X,Y)=1$. By the invariance of the Killing form $B$
we get
$$
1 = B(X,Y) = B(\Ad(t)X,\Ad(t)Y) = cdB(X,Y) =cd.
$$
So on each pair of root spaces $\Ad(t)$ has determinant one hence
also on $\p$.

Replacing $G$ by a double cover if necessary, which doesn't effect
the claim of the lemma, we may assume that $\ph$ lifts to the spin
group $\Spin(\p)$. Let $\p=\p^+\oplus\p^-$ be the decomposition
according to an ordering of $\phi(\t ,\g)$. This decomposition is
a polarization of the quadratic space $\p$ and hence the spin
group acts on $S^+=\wedge^{even}\p^+$ and $S^-=\wedge^{odd}\p^+$
in a way that the virtual module $(S^+-S^-)\otimes (S^+-S^-)$
becomes isomorphic to $W$. For $K$ connected the claim now follows
from \cite{AtSch} (4.5). An inspection shows however that the
proof of (4.5) in \cite{AtSch}, which is located in the appendix
(A.12), already applies when we only assume that the homomorphism
$\ph$ factors over the spin group. \qed

We continue the proof of the theorem. Let $\hat{T}$ denote the
set of all unitary characters of $T$. Any regular element
$\hat{t}\in \hat{T}$ gives rise to a discrete series
representation $(\omega ,V_\omega)$ of $G$. Let
$\Theta_{\hat{t}}=\Theta_\omega$ be its character which, due to
Harish-Chandra, is known to be a function on $G$. Harish-Chandra's
construction gives a bijection between the set of discrete series
representations of $G$ and the set of $W(G,T)=W(K,T)$-orbits of
regular characters of $T$.

Let $\Phi^+$ denote the set of positive roots of $(\g,\t)$ and let
$\Phi^+_c$ ,$\Phi^+_n$ denote the subsets of compact and
noncompact positive roots. For each root $\alpha$ let $t\mapsto
t^\alpha$ denote the corresponding character
 on $T$.
Define
\begin{eqnarray*}
'\Delta_c(t) &=& \prod_{\alpha\in\Phi^+_c}(1-t^{-\alpha})\\
'\Delta_n &=& \prod_{\alpha\in\Phi^+_n}(1-t^{-\alpha})
\end{eqnarray*}
and $'\Delta = '\Delta_c '\Delta_n$. If $\hat{t}\in \hat{T}$ is
singular, Harish-Chandra has also constructed an invariant
distribution $\Theta_{\hat{t}}$ which is a virtual character on
$G$. For $\hat{t}$ singular let $W(\hat{t})\subset W(\g,\t)$ be
the isotropy group. One has $\Theta_{\hat{t}} = \sum_{w\in
W(\hat{t})} \epsilon(w)\Theta'_{w,\hat{t}}$ with
$\Theta'_{w,\hat{t}}$ the character of an induced representation
acting on some Hilbert space $V_{w,\hat{t}}$ and $\epsilon(w)\in
\{\pm 1\}$. Let $\CE_2(G)$ denote the set of discrete series
representations of $G$ and $\CE_2^s(G)$ the set of $W(G,T)$-orbits
of singular characters.

By Theorem \ref{inv-orb-int} the current theorem will follow from the

\begin{lemma}\label{trace-tau}
For $t\in T$ regular we have
$$
\tr \tau(t)\=
\rez{|W(G,T)|}\sum_{\hat{t}\in\hat{T}}\Theta_{\hat{t}}(f_\tau)
\Theta_{\hat{t}}(t).
$$
\end{lemma}

\prf Let $\ga$ denote the virtual $K$-representation on
$(\wedge^{even}\p -\wedge^{odd}\p)\otimes V_\tau$. Harish-Chandra
has shown (\cite{HC-S} Theorem 12) that for any
$\hat{t}\in\hat{T}$ there is an irreducible unitary representation
$\pi_{\hat{t}}^0$ such that $\Theta_{\hat{t}}$ coincides up to
sign with the character of $\pi_{\hat{t}}^0$ on the set of
elliptic elements of $G$ and $\pi_{\hat{t}}^0=\pi_{\hat{t}'}^0$ if
and only if there is a $w\in W(G,T)=W(K,T)$ such that
$\hat{t}'=w\hat{t}$.

Further (\cite{HC-S}, Theorem 14) Harish-Chandra has shown that
the family
$$
\left( \frac{'\lap(t)\Theta_{\hat{t}}(t)}
            {\sqrt{|W(G,T)|}}\right)_{\hat{t}\in \hat{T}/W(G,T)}
$$
forms an orthonormal basis of $L^2(T)$. Here we identify
$\hat{T}/W(G,T)$ to a set of representatives in $\hat{T}$ to make
$\Theta_{\hat{t}}$ well defined.

Consider the function $g(t)=\frac{\tr\ga(t) \
'\lap_c(t)}{\overline{\ '\lap_n(t)}} = \tr\tau(t) \ '\lap(t)$. Its
coefficients with respect to the above orthonormal basis are
\begin{eqnarray*}
\langle g, \frac{\ '\lap\Theta_{\hat{t}}}
            {\sqrt{|W(G,T)|}}\rangle &=&
            \rez{\sqrt{|W(G,T)|}}\int_T\tr\ga(t) |
            '\lap_c(t)|^2\overline{\Theta_{\hat{t}}(t)} dt\\
 &=& \sqrt{|W(G,T)|}\int_K \tr\ga(k)\overline{\Theta_{\hat{t}}(k)}
 dk
\end{eqnarray*}
where we have used the Weyl integration formula for the group $K$
and the fact that $W(G,T)=W(K,T)$. Next by Lemma
\ref{fin-lin-comb} this equals
$$
\sqrt{|W(G,T)|}
\dim((\wedge^{even}\p-\wedge^{odd}\p)\otimes\breve{\tau}\otimes\pi_{\hat{t}}^0)^K
\= \sqrt{|W(G,T)|}\Theta_{\hat{t}}(f_\tau).
$$
Hence
\begin{eqnarray*}
g(t) &=& \tr\tau(t) \ '\lap(t)\\
    &=& \sum_{\hat{t}\in \hat{T}/W(G,T)}\Theta_{\hat{t}}(f_\tau)
    \ '\lap(t) \Theta_{\hat{t}}(t)\\
    &=& \rez{|W(G,T)|}\sum_{\hat{t}\in \hat{T}}\Theta_{\hat{t}}(f_\tau)
    \ '\lap(t) \Theta_{\hat{t}}(t).
\end{eqnarray*}
The lemma and the theorem are proven. \qed

\begin{corollary}
If $\tilde{g}\in\tilde{G}$ is semisimple and not elliptic then
$\CO_{\tilde{g}}(g_\tau)=0$. If $\tilde{g}$ is elliptic regular
then
$$
\CO_{\tilde{g}}(g_\tau)\=\frac{\tr\tau(\tilde{g})} {\tr(\tilde{g}
| S^+-S^-)}.
$$
\end{corollary}

\prf Same as for the last proposition with $g_\tau$ replacing
$f_\tau$.
\qed

\begin{proposition}
Assume that $\tau$ extends to a representation of the group $G$
on the same space. For the function $f_{\tau}$ we have for any
$\pi \in \hat{G}$:
$$
\tr \ \pi(f_{{\tau}}) \= \sum_{p=0}^{\dim \ \g /\k}(-1)^p \dim \
{\rm Ext}_{(\g ,K)}^p (V_{{\tau}} ,V_\pi),
$$
i.e., $f_{{\tau}}$ gives the Euler-Poincar\'e numbers of the $(\g
,K)$-modules $(V_{{\tau}} ,V_\pi)$, this justifies the name
Euler-Poincar\'e function.
\end{proposition}

\prf By definition it is clear that
$$
\tr \ \pi (f_{{\tau}}) \= \sum_{p=0}^{\dim \ \p} (-1)^p \dim \
H^p(\g ,K,V_{\breve{\tau}} \otimes V_\pi).
$$
The claim now follows from \cite{BorWall}, p. 52. \qed

\section{The Selberg trace formula}
In this section we will fix the basic notation and set up
the trace formula. For compactly supported functions this
formula is easily deduced. In the sequel we however will
need it for functions with noncompact support and
therefore will have to show more general versions of the
trace formula.

Let $G$ denote a real reductive group.

\subsection{The trace formula}
Let $\Ga\subset G$ be a discrete subgroup such that
the quotient manifold $\Ga\bs G$ is compact. We say that $\Ga$ is
{\it cocompact} in $G$. Examples are given by
nonisotropic arithmetic groups \cite{marg}. Since $G$ is
unimodular the Haar measure on $G$ induces a $G$-invariant measure
on $\Ga\bs G$, so we can form the Hilbert space $L^2(\Ga\bs G)$
 of square integrable measurable functions
on $\Ga\bs G$ modulo null functions. More generally, let
$(\omega,V_\omega)$ be a finite
dimensional unitary representation of $\Ga$ and let $L^2(\Ga\bs
G,\omega)$ be the Hilbert space of
all measurable functions $f : G\ra V_\omega$ such that $f(\ga
x)=\omega(\ga)f(x)$ for all $\ga\in\Ga$ and all $x\in G$ and such
that
$$
\int_{\Ga\bs G} \norm{f(x)}^2 dx\ <\ \infty
$$
modulo null functions. The scalar product of $f,g\in L^2(\Ga\bs
G,\omega)$ is
$$
(f,g) \= \int_{\Ga\bs G} \langle f(x),g(x)\rangle dx,
$$
where $\langle .,.\rangle$ is the scalar product on $V_\omega$.
Let $C^\infty(\Ga\bs G,\omega)$
be the subspace consisting of smooth functions.

The group $G$ acts unitarily on $L^2(\Ga\bs G,\omega)$ by
$$
R(y)\ph(x)\= \ph(xy)
$$
for $x,y\in G$ and $\ph\in L^2(\Ga\bs G,\omega)$. Let $\pi$ be
any unitary representation of $G$, any $f\in L^1(G)$ defines a
bounded operator
$$
\pi(f) \= \int_G f(x) \pi(x) dx
$$
on the space of $\pi$. We apply this to the case
$\pi=R$. Let $C_c^\infty(G)$ be the space of all smooth functions
of compact support on $G$. Let $f\in C_c^\infty(G)$, $\ph\in
L^2(\Ga\bs G,\omega)$. Fix a fundamental domain $\CF\subset G$
for the $\Ga$-action on $G$ and compute formally at first:
\begin{eqnarray*}
R(f)\ph(x) &=& \int_G f(y) \ph(xy) dy\\
    &=& \int_G f(x^{-1}y) \ph(y) dy\\
    &=& \sum_{\ga\in\Ga} \int_{\ga\CF}f(x^{-1}y)\ph(y)dy\\
    &=& \sum_{\ga\in\Ga}\int_\CF f(x^{-1}\ga y) \omega(\ga)\ph(y) dy\\
    &=& \int_{\Ga\bs G} k_f(x,y)\ph(y) dy,
\end{eqnarray*}
where $k_f(x,y) =\sum_{\ga\in\Ga}f(x^{-1}\ga y)\omega(\ga)$.
Since $f$ has compact support the latter sum is locally finite
and therefore defines a smooth Schwartz kernel on the compact
manifold $\Ga\bs G$. This implies that the operator $R(f)$ is a
smoothing operator and hence of trace class.

Since the convolution algebra $C_c^\infty$ contains an approximate
identity we can infer that the unitary $G$-representation $R$ on
$L^2(\Ga\bs G,\omega)$ decomposes into a direct sum of
irreducibles, i.e.,
$$
L^2(\Ga\bs G,\omega) \ \cong\ \bigoplus_{\pi\in\hat{G}} N_{\Ga
,\omega}(\pi)\pi
$$
with finite multiplicities $N_{\Ga ,\omega}(\pi)\in\N_0$.
 Moreover, the trace of $R(f)$ is
given by the integral over the diagonal, so
$$
\tr R(f) \= \int_{\Ga\bs G} \tr k_f(x,x) dx,
$$
where the trace on the right hand side is the trace in
$\End(V_\omega)$. We plug in the sum for $k_f$ and rearrange that
sum in that we first sum over all conjugacy classes in the group
$\Ga$. We write $\Ga_\ga$ and $G_\ga$ for the centralizers of
$\ga\in\Ga$ in $\Ga$ and in $G$ resp.
\begin{eqnarray*}
\tr R(f) &=& \sum_{\ga\in\Ga} \tr\omega(\ga) \int_\CF f(x^{-1}\ga x) dx\\
    &=& \sum_{[\ga]}\sum_{\sigma\in\Ga /\Ga_\ga} \tr\omega(\ga)
        \int_\CF f((\sigma x)^{-1}\ga \sigma x) dx\\
    &=& \sum_{[\ga]}\tr\omega(\ga)
        \int_{\Ga_\ga\bs G} f(x^{-1}\ga x) dx\\
    &=& \sum_{[\ga]}\tr\omega(\ga)\vol(\Ga_\ga\bs G_\ga) \CO_\ga(f),
\end{eqnarray*}
where $\CO_\ga(f)$ is the orbital integral. We have proved the following
proposition.

\begin{proposition}
For $f\in C_c^\infty(G)$ we have $$
\sum_{\pi\in\hat{G}} N_{\Ga,\omega}(\pi)\tr\pi(f) \= \sum_{[\ga]}
\tr\omega(\ga) \vol(\Ga_\ga\bs G)\CO_\ga(f),
$$
where all sums and integrals converge absolutely.\qed
\end{proposition}

All this is classical and may be found at various places. For our
applications we will need to extend the range of functions $f$ to
be put into the trace formula. For this sake we prove the proposition below.

Let $k,l\in\N$ and define $L^1_{k,l}(G)$ to be the set of all functions $f$ on $G$ which are $\max{k,l}$-times continuously differentiable and satisfy
\begin{itemize}
\item $|Df|$ is integrable on $G$ for every $D\in U(\g)$ of degree $\le k$, and
\item $|Df|$ is bounded on $G$ for every $D\in U(\g)$ of degree $\le l$.
\end{itemize}
Every $D\in U(\g)$  of degree $\le k$ induces a seminorm on $L^1_{k,l}(G)$ by
$$
\sigma_D(f)\= \int_G|Df(x)|\, dx.
$$
Further, every $D\in U(\g)$  of degree $\le l$ induces a seminorm on $L^1_{k,l}(G)$ by
$$
s_D(f) \= \sup_{x\in G}|Df(x)|.
$$
where $\norm{\cdot}_1$ is the $L^1$-norm.
We equip $L^1_{k,l}(G)$ with the topology given by these seminorms.

\begin{proposition}\label{convergence-of-tf}
Suppose $f\in L_{2N,1}^1(G)$ with $N>\frac{\dim G}{2}$.
 Then the trace formula is valid for
$f$ and either side of the trace formula defines a continuous linear functional on $L_{2N,1}^1(G)$.
\end{proposition}

\prf  We consider $U(\g)$ as the algebra of all left invariant
differential operators on $G$. Choose a left invariant Riemannian
metric on $G$ and let $\lap$ denote the corresponding Laplace
operator. Then $\lap\in U(\g)$ and thus it makes sense to write
$R(\lap)$, which is an elliptic differential operator of order
$2$ on the compact manifold $\Ga\bs G$, essentially selfadjoint
and non-negative. The theory of pseudodifferential operators
implies that $R(\lap +1)^{-N}$ has a $C^1$-kernel and thus is of
trace class. Let $g=(\lap +1)^Nf$ then $g\in L^1(G)$, so $R(g)$
is defined and gives a continuous linear operator on the Hilbert
space $L^2(\Ga \bs G,\ph)$. We infer that $R(f)=R(\lap
+1)^{-N}R(g)$ is of trace class.

Let $\chi :[0,\infty [\ra [0,1]$ be a monotonic $C^{2N}$-function
with compact support, $\chi\equiv 1$ on $[0,1]$ and
$|\chi^{(k)}(t)|\le 1$ for $k=1,\dots ,2N$. Let
$h_n(x):=\chi(\frac{dist(x,e)}{n})$ for $x\in G$, $n\in \N$. Then
$|Dh_n(x)|\le \frac{C_D}{n}$ for any $D\in\g U(\g)$. Let
$f_n=h_nf$ then $f_n\ra f$locally uniformly.

{\it Claim.} For the $L^1$-norm on $G$ we have
$$
\parallel Df_n -Df \parallel_1 \ra 0
$$
as $n\ra \infty$ for any $D\in U(\g)$ of degree $\le 2N$.

Proof of the claim: By the Poincar\'e-Birkhoff-Witt theorem
$D(f_n)=D(h_nf)$ is a sum of expressions of the type
$D_1(h_n)D_2(f)$ and $D_1$ can be chosen to be the identity
operator or in $\g U(\g)$. The first case gives the summand
$h_nD(f)$ and it is clear that $\norm{Df-h_nDf}_1$ tends to zero
as $n$ tends to infinity. For the rest assume $D_1\in \g U(\g)$.
Then $|D_1(h_n)(x)|\le \frac{c}{n}$ hence
$\norm{D_1(h_n)D_2(f)}_1$ tends to zero because $D_2(f)$ is in
$L^1(G)$. The claim follows.

To prove the lemma we estimate the operator norm as
$$
\norm{R((\lap +1)^Nf_n)-R((\lap +1)^Nf)}\ \le\ \norm{(\lap
+1)^Nf_n-(\lap +1)^Nf}_1
$$
the latter tends to zero according to the claim proven above.
Denoting the trace norm by $\norm{.}_{tr}$ we infer
$$
\norm{R(f_n)-R(f)}_{tr}
$$ $$
= \norm{R(\lap +1)^{-N}(R((\lap +1)^Nf_n)-R((\lap +1)^Nf))}_{tr}
$$ $$
\le \norm{R(\lap +1)^{-N}}_{tr}\norm{R((\lap +1)^Nf_n)-R((\lap
+1)^Nf)}
$$
which tends to zero. Therefore $\tr R(f_n)$ tends to $\tr R(f)$
as $n\ra \infty$. It follows
\begin{eqnarray*}
\sum_{\pi\in\hat{G}} N_{\Ga ,\omega}(\pi)\tr\pi (f) &=& \tr R(f)\\
    &=& \lim_n \tr R(f_n)\\
    &=& \lim_n \sum_{[\ga]} \tr\omega(\ga)\vol(\Ga_\ga\bs G_\ga) \CO_\ga(f_n).
\end{eqnarray*}
Now suppose as an additional condition that $f\ge 0$ and $\omega =1$.
Then we are allowed to interchange the limit and the sum by monotone convergence and thus in this case
$$
\sum_{\pi\in\hat{G}} N_{\Ga ,1}(\pi)\tr\pi (f)\=
\sum_{[\ga]} \vol(\Ga_\ga\bs G_\ga) \CO_\ga(f).
$$
In particular, the right hand side is finite.
For general $f$ and $\omega$, we use the boundedness of $\tr\omega$ and we can justify the same interchange by dominated convergence if we show that there is $\tilde f\in L_{2N,1}^1(G)$ with $\tilde f\ge |f|$, because then the trace formula is valid for $\tilde f$ and thus 
$$
\sum_{[\ga]} \vol(\Ga_\ga\bs G_\ga) \CO_\ga(\tilde f)\ <\ \infty.
$$
It remains to show the existence of $\tilde f$.
For this let $U,V$ be small neighbourhoods of the unit in $G$ with $\overline U\subset V$.
Let $\al\in C^\infty(G)$ with values in $[0,1]$, support in $V$ and such that $\al\equiv 1$ on $U$.
Since $Xf$ is bounded for every $X\in\g$ it follows that there is $C>0$ such that for small $t>0$ and every $x\in G$ one has
\begin{eqnarray*}
||f(x\exp(tX))|-|f(x)|| &\le& |f(x\exp(tX))-f(x)|\\
&\le& Ct.
\end{eqnarray*}
One has
\begin{eqnarray*}
|f|*\al(x) &=& \int_G|f(y)|\al(y^{-1}x)\, dy\\
&=& \int_G|f(xy^{-1}|\al(y)\, dy\\
&\ge& \int_U|f(xy^{-1})|\, dy.
\end{eqnarray*}
It follows that there exists $\eps >0$ with
$$
|f| *\al(x)\ \ge\ \eps|f(x)|
$$
for every $x\in G$.
Set $\tilde f=\frac 1\eps |f|*\al$, then $\tilde f$ lies in $L_{2N,1}^1(G)$ and $\tilde f\ge |f|$.
By the above the trace formula is valid for $f$.
Finally, the norm-estimates also imply the claimed continuity of the linear functional $f\mapsto \tr(R(f))$ on $L_{2N}^1(G)$.
\qed

\section{The Lefschetz formula}

In this section $G$ will be a connected semisimple Lie
group with finite center.

\subsection{Euler characteristics}
Let $L$ be a real reductive group and suppose there is a finite
subgroup $E$ of the center of $L$ and a reductive and Zariski-connected linear group
 $\CL$ over $\R$ such that $L/E$ is
isomorphic to a subgroup of $\CL(\R)$ of finite index. Note that
these conditions are satisfied whenever $L$ is a Levi component of
a connected semisimple group $G$ with finite center. Let $K_L$ be
a maximal compact subgroup of $L$ and let $\Ga$ be a cocompact
discrete subgroup of $L$. Fix a nondegenerate invariant bilinear
form on the Lie algebra $\l_0$ of $L$ such that $B$ is negative
definite on the Lie algebra of $K_L$ and positive definite on its
orthocomplement. Let $\theta$ be the Cartan involution fixing
$K_L$ pointwise then the form $-B(X,\theta(Y))$ is positive
definite and thus defines a left invariant metric on $L$. For any
closed subgroup $Q$ we get a left invariant metric on $Q$. The
volume element of that metric gives a Haar measure, called the
{\it standard volume}  with respect to $B$
on $Q$.

Let $\Ga\subset G$ denote a cocompact discrete subgroup. If $\Ga$ is torsion free, it acts fixed point free
on the contractible space $X$ and hence $\Ga$ is the fundamental
group of the Riemannian manifold
$$
X_{\Ga} = \Ga \bs X = \Ga \bs G/K
$$
 it follows that we have a
canonical bijection of the homotopy classes of loops:
$$
[S^1 : X_{\Ga} ] \ra \Ga / {\rm conjugacy}.
$$
For a given class $[\ga]$ let $X_\ga$ denote the
union of all closed geodesics in the corresponding class in $[S^1
: X_\Ga ]$. Then $X_\ga$ is a smooth submanifold of $X_{\Ga_H}$
\cite{DKV}, indeed, it follows that
$$
X_\ga\ \cong\ \Ga_\ga \bs G_\ga/K_\ga,
$$
where $G_\ga$ and $\Ga_\ga$ are the centralizers of $\ga$ in $G$
and $\Ga$ and $K_\ga$ is a maximal compact subgroup of $G_\ga$.
Further all closed geodesics in the class $[\ga]$ have the same
length $l_\ga$.

Let $r\in\N_0$. If $\Ga$ is torsion-free, we define the $r$-th Euler characteristic of $\Ga$ by
$$
\chi_r(\Ga)\=\chi_r(X_\Ga)\=\sum_{j=0}^{\dim X_\Ga} (-1)^{j+r} \left(\begin{array}{c}j \\r\end{array}\right) b^j(X_\Ga),
$$
where $b^j(X_\Ga)$ is the $j$-th Betti number of $X_\Ga$.
We want to extend this notion to groups $\Ga$ which are not necessarily torsion-free.

Let now $H$ be a $\theta$-stable Cartan subgroup, then $H=AB$,
where $A$ is a connected split torus and $B\subset K_L$ is a
Cartan of $K_L$. On the space $L/H$ we have a pseudo-Riemannian
structure given by the form $B$. The Gauss-Bonnet construction
(\cite{dieu} sect. 24 or see below) generalizes to
Pseudo-Riemannian structures to give an Euler-Poincar\'e measure
$\eta$ on $L/H$. Define a (signed) Haar-measure by
$$
\mu_{EP} \= \eta\otimes ({\rm normalized\ Haar\ measure\ on}\ H).
$$
Let $W=W(L,H)$ denote the Weyl group and let $W_\C=W(L_\C,H_\C)$
be the Weyl group of the complexifications. Let the {\it generic
Euler characteristic} be
defined by
$$
\chi_{gen}(\Ga \bs L/K_L)\= \frac{\mu_{EP}(\Ga\bs L)}{|W|}.
$$
Write $X_\Ga = \Ga\bs L/K$.

\begin{lemma}
If $H$ is compact and $\Ga$ is torsion-free, then the generic
Euler characteristic equals the ordinary Euler characteristic,
i.e., $\chi_{gen}(X_\Ga)=\chi(X_\Ga)$.
\end{lemma}

\prf In this case we have $H=B$. The Gauss-Bonnet Theorem tells us
$$
\eta(\Ga\bs L/B)\= \chi(\Ga\bs L/B).
$$
Now $\Ga\bs L/B\ra \Ga\bs L/K_L$ is a fiber bundle with fiber
$K_L/B$, therefore
 we get
$$
\chi(\Ga\bs L/B)\= \chi(\Ga\bs L/K_L)\chi(K_L/B).
$$
Finally the Hopf-Samelson formula says $\chi(K_L/B)=|W|$. \qed

For the next proposition assume that $A$ is central in $L$, $L=AL_1$ where $L_1$ has compact
center. Let $C$ denote the center of $L$, then
$A\subset C$ and $C=AB_C$, where $B_C=B\cap C$. Let $L'$ be the
derived group of $L$ and let $\Ga'=L'\cap \Ga C$ and
$\Ga_C=\Ga\cap C$ then by Lemma 3.3 in \cite{Wolf} we infer that
$\Ga_C$ is a cocompact subgroup of $C$ and $\Ga'$ is
a cocompact discrete subgroup of $L'$. 
Let $\Ga_A=A\cap \Ga_C
B_C$ the projection of $\Ga_C$ to $A$. Then $\Ga_A$ is discrete
and cocompact in $A$.

\begin{proposition}\label{4.4}
Assume $\Ga$ is  torsion-free and $A$ is central in $L$ of dimension $r$, and $\Ga'\subset\Ga$ is a subgroup of finite index. Then the group
$A/\Ga_A$ acts freely on $X_\Ga$ and $\chi_{gen}(X_\Ga)=\chi_r(\Ga)\vol(A/\Ga_A)$.
It follows that
$$
\chi_r(\Ga)\=\chi_r(\Ga')\,\frac{[\Ga_A:\Ga_A']}{[\Ga:\Ga']}.
$$
\end{proposition}

\prf The group $A_\Ga=A/\Ga_A$ acts on $\Ga \bs L/B$ by
multiplication from the right. We claim that this action is free,
i.e., that it defines a fiber bundle
$$
A_\Ga\ra\Ga\bs L/B\ra \Ga\bs L/H.
$$
To see this let $\Ga xaB=\Ga xB$ for some $a\in A$ and $x\in L$,
then $a=x^{-1}\ga xb$ for some $\ga\in\Ga$ and $b\in B$. Writing
$\ga$ as $\ga' a_\ga b_\ga$ with $\ga'\in\Ga'$ and $a_\ga\in A$
and $b_\ga\in B$ we conclude that $a_\ga\in\Ga_A$ and $a=a_\ga$, whence the claim.

In the same way we see that we get a fiber bundle
$$
A_\Ga\ra\Ga\bs L/K_L\ra A\Ga\bs L/K_L.
$$
We now apply the Gauss-Bonnet theorem to conclude
\begin{eqnarray*}
\chi_{gen}(X_\Ga) &=& \vol(A/\Ga_A)\frac{\eta(\Ga\bs L/H)}{|W|}\\
    &=& \vol(A/\Ga_A)\frac{\chi(\Ga\bs L/H)}{\chi(K_L/B)}\\
    &=& \vol(A/\Ga_A)\frac{\chi(A\Ga\bs L/B)}{\chi(K_L/B)}\\
    &=& \vol(A/\Ga_A)\chi(A\Ga\bs L/K_L).
\end{eqnarray*}
It remains to show that $\chi(A\bs X_\Ga)=\chi_r(\Ga)$.
For this let $\a_0$ be the real Lie algebra of $A$ and let $\l_0^{\rm der}$ be the Lie algebra of the derived group $L^{\rm der}$.
The Lie algebra $\l_0$ of $L$ can be written as
$$
\l_0\=\a_0\oplus \l_0^{\rm der}\oplus\z_0,
$$
where $\z_0$ is central in $\l_0$.
Let $X_1,\dots,X_r$ be a basis of $\a_0$.
We consider $X_j$ as a vector filed on $\Ga\bs L/K_L$ by means of the left translation.
Let $\omega_1,\dots,\omega_r$ be the dual basis of $\a_0^*$.
Via the above decomposition we can view each $\omega_j$ as an element of $l_0^*$, thus as a $1$-form on $L$ which is bi-invariant and induces a $1$-form on $\Ga\bs L/K_L$.
Since $A\cong\R^r$ and the $\omega_j$ are the differential forms given by a set of co-ordinates, the forms $\omega_1,\dots,\omega_r$ are all closed.
The group $A_\Ga=A/\Ga_A$ is connected and compact, therefore the cohomology of the deRham complex of $\Ga\bs L/K_L$ coincides with the cohomology of the subcomplex of $A_\Ga$-invariants $\Omega(X_\Ga)^{A_\Ga}$.
Using local triviality of the bundles one sees that
$$
\Omega(X_\Ga)^{A_\Ga}\=\bigoplus_{I\subset\{ 1,\dots,r\}}\pi^*\Omega(A\bs X_\Ga)\wedge\omega_I,
$$
where $\pi$ is the projection $X_\Ga =\Ga\bs L/K_L\to A\Ga\bs L/K_L=A\bs X_\Ga$ and 
$$
\omega_{\{i_1,\dots,i_k\}}=\omega_{i_1}\wedge\dots\wedge\omega_{i_k}
$$ 
for $i_1<i_2<\dots <i_k$.
Since the $\omega_j$ are closed, it follows for the real valued cohomology that
$$
H^\bullet(X_\Ga)\ \cong\ \bigoplus _{I\subset\{ 1,\dots,r\}}H^{\bullet- |I|}(A\bs X_\Ga).
$$
So we compute
\begin{eqnarray*}
\chi_r(\Ga) &=& \sum_{j=r}^{\dim X_\Ga} (-1)^{j+r}\left(\begin{array}{c}j \\r\end{array}\right) b^j(X_\Ga)\\
&=& \sum_{j=r}^{\dim X_\Ga} (-1)^{j+r}\left(\begin{array}{c}j \\r\end{array}\right) \sum_{I\subset\{1,\dots,r\}} b^{j-|I|}(A\bs X_\Ga)\\
&=& \sum_{j=r}^{\dim X_\Ga} (-1)^{j+r}\left(\begin{array}{c}j \\r\end{array}\right) \sum_{k=0}^r\left(\begin{array}{c}r \\k\end{array}\right)b^{j-k}(A\bs X_\Ga)\\
&=& \sum_{p=0}^{\dim A\bs X_\Ga} b^p(A\bs X_\Ga)\sum_{j=r}^{\dim X_\Ga}(-1)^{j+r}\left(\begin{array}{c}j \\r\end{array}\right)\left(\begin{array}{c}r \\j-p\end{array}\right)\\
&=&\chi(A\bs X_\Ga).
\end{eqnarray*}
The last step uses the combinatorial identity
$$
\sum_{j=r}^{r+p}(-1)^{j+r}\left(\begin{array}{c}j \\r\end{array}\right)\left(\begin{array}{c}r \\j-p\end{array}\right)\= (-1)^p.
$$
\qed

Define the $r$-th Euler-number of $\Ga$ by
$$
\chi_r(\Ga)\=\chi_r(\Ga')\,\frac{[\Ga_A:\Ga_A']}{[\Ga:\Ga']},
$$
where $\Ga'\subset\Ga$ is a torsion-free subgroup of finite index (which always exists by Selberg's Lemma \cite{Sel2}).
Proposition \ref{4.4} shows that $\chi_r(\Ga)$ does not depend on the choice of $\Ga'$.
Further, this definition allows us to extend Proposition \ref{4.4} to arbitrary cocompact $\Ga$:
$$
\chi_{\rm gen}(X_\Ga)\=\chi_r(\Ga)\,\vol(A/\Ga_A).
$$

We will compute $\chi_{gen}(X_\Ga)$ in terms of root systems. Let
$\Phi$ denote the root system of $(\l,\h)$, where $\l$ and $\h$
are the complexified Lie algebras of $L$ and $H$. Let $\Phi_n$ be
the set of noncompact imaginary roots and choose a set $\Phi^+$ of
positive roots such that for $\alpha\in\Phi^+$ nonimaginary we
have that $\alpha^c\in\Phi^+$. Let $\nu=\dim L/K_L-\rank L/K_L$,
let $\rho$ denote the half of the sum of all positive roots. For
any compact subgroup $U$ of $L$ let $v(U)$ denote the standard
volume.

\begin{theorem} \label{3.3}
The generic Euler number satisfies
\begin{eqnarray*}
\chi_{gen}(X_\Ga)&=& \frac{(-1)^{|\Phi_n^+|}|W_\C|
                \prod_{\alpha\in\Phi^+}(\rho,\alpha)v(K_L)}
            {(2\pi)^{|\Phi^+|}|W|2^{\nu/2}v(B)}
\vol(\Ga\bs L)\\
    &=& c_L^{-1}|W_\C|\prod_{\alpha\in\Phi^+}(\rho,\alpha)\vol(\Ga\bs L),
\end{eqnarray*}
where $c_L$ is Harish-Chandra's constant, i.e.,
$$
c_L\= (-1)^{|\Phi_n^+|}
(2\pi)^{|\Phi^+|}2^{\nu/2}\frac{v(B)}{v(K_L)}|W|.
$$
So, especially in the case when $A$ is
central we get
$$
\chi_r(\Ga)\= \frac{|W_\C|\prod_{\alpha\in\Phi^+}
(\rho,\alpha)}{c_L\vol(A/\Ga_A)}\vol(\Ga\bs L).
$$
\end{theorem}

\prf On the manifold $L/H$ the form $B$ gives the structure of a
pseudo-Riemannian manifold. Let $P$ denote the corresponding
$SO_{p,q}$ fiber bundle, where $(p,q)$ is the signature of $B$ on
$\l_0/\h_0$. Let $\ph :H\ra SO_{p,q}$ denote the homomorphism
induced by the adjoint representation. We the have $P=L\times_\ph
SO_{p,q}$. A connection on $P$ is given by the $L$-invariant
connection $1$-form
$$
\omega \left( A+\sum_\alpha c_\alpha X_\alpha \right) \= A,\ \ \
\ A\in so(p,q),
$$
where we have used $B_eP\cong
so(p,q)\oplus\left((\oplus_\alpha\g_\alpha)\cap\g_0\right)$. By an
inspection in local charts one finds the following formula for the
$L$-invariant $2$-form $d\omega$:
$$
d\omega(e)\left( A+\sum_\alpha c_\alpha X_\alpha,A'+ \sum_\alpha
c_\alpha' X_\alpha\right) = [A,A'] -\sum_\alpha c_\alpha
c_\alpha' \ph_*H_\alpha.
$$
Let $\Omega=d\omega$ and let $Pf$ be the Pfaffian as in
\cite{dieu} 24.46.10. Let $J$ be the diagonal matrix having the diagonal entries
$(1,\dots,1,-1,\dots,-1)$ with $p$ ones and $q$ minus
ones and let $P(X)=Pf(JX)$. Write this as
$$
P(X)\= (2\pi)^{-m}\sum_{h,k}\chi_{h,k}X_{h_1,k_1}\dots
X_{h_m,k_m},
$$
where $m-\rez{2}(p+q)\in\N$. Let
$$
F(\Omega)\= (2\pi)^{-m} \sum_{h,k}
\chi_{h,k}\Omega_{h_1,k_1}\wedge\dots\wedge\Omega_{h_m,k_m}.
$$
There is a unique $G$-invariant form $F_B(\Omega)$ on $L/H$ which
at the origin is the pullback of $F(\Omega)$ with respect to a
section $s: L/H\ra P$. Writing $\Omega =\Omega_1 -\Omega_2$ with
$$
\Omega_1\left( A+\sum_\alpha c_\alpha X_\alpha\ ,\ A'+ \sum_\alpha
c_\alpha' X_\alpha\right) = [A,A']
$$
we get
$$
F_B(\Omega)\= F_B(-\Omega_2)\= (-1)^m F_B(\Omega_2).
$$
We call this form $\eta$. On the space $\C X_\alpha +\C
X_{-\alpha}$ we have that
$$
\ph_*Y=-\alpha(Y)i\matrix{{}}{1}{1}{{}}
$$
with respect to the
basis $(X_\alpha,X_{-\alpha})$. Let $\alpha_1,\dots,\alpha_m$ be
an enumeration of $\Phi^+$, then
$$
\eta\= (2\pi)^{-m} \sum_{\sigma\in Per(\Phi)}\sum_{h,k}\chi_{h,k}
(\ph_*H_{\sigma(\alpha_1)})_{h_1,k_1}\dots
(\ph_*H_{\sigma(\alpha_m)})_{h_m,k_m} \tilde{\omega}_\sigma,
$$
where
$$
\tilde{\omega}_\sigma\= dX_{\sigma(\alpha_1)}\wedge
dX_{-\sigma(\alpha_1)} \wedge\dots\wedge
dX_{\sigma(\alpha_m)}\wedge dX_{-\sigma(\alpha_m)}.
$$
We end up with
\begin{eqnarray*}
\eta &=& (2\pi)^m(-1)^m \sum_{{\alpha\in\Phi^+},\
{\sigma\in\Per(\Phi^+)}}
(\alpha,\sigma(\alpha))\tilde{\omega}_\sigma\\
    &=& (2\pi)^m(-1)^m |W_\C|\prod_{\alpha\in\Phi^+} (\rho,\alpha)
\tilde{\omega}_\sigma.
\end{eqnarray*}
Comparing $\eta$ with the standard measure given by the form $B$
gives the claim.
\qed

\subsection{Setting up the formula}\label{4.1}
Let $K\subset G$ be a maximal compact subgroup with Cartan
involution $\theta$. Let $X=G/K$ denote the symmetric space
attached to $G$. Let $H\subset G$ be a Cartan subgroup. Modulo
conjugation we may assume that $H$ is stable under $\theta$. Then
$H=AB$, where $A$ is a connected split torus and $B$ is a subgroup
of $K$. The double use of the letter $B$ for the group and a
bilinear form on the Lie algebra will not cause any confusion. Fix
a parabolic subgroup $P$ of $G$ with split component $A$. Then
$P=MAN$ where $N$ is the unipotent radical of $P$ and $M$ is
reductive with compact center and finite component group. The
choice of the parabolic $P$ amounts to the same as a choice of a
set of positive roots $\Phi^+=\Phi^+ (\g ,\a)$ in the root system
$\Phi(\g,\a)$ such that for the Lie algebra $\n={\rm Lie}_\C(N)$ we have
$\n=\bigoplus_{\alpha\in\Phi^+}\g_\alpha$. Let
$\bar{\n}=\bigoplus_{\alpha\in\Phi^+}\g_{-\alpha}$, $\bar\n_0=\bar\n\cap\g_0$ and
$\bar{N}=\exp(\bar{\n}_0)$ then $\bar{P}=MA\bar{N}$ is the parabolic
{\it opposite} to $P$. The root space
decomposition then writes as
$\g=\a\oplus\m\oplus\n\oplus\bar{\n}$. Let $\rho_P$ be the half of
the sum of the positive roots, each weighted with its
multiplicity, i.e., for $a\in A$ we have $a^{2\rho_P}
=\det(a|\n)$. Let $A^-\subset A$ denote the negative
Weyl chamber corresponding to that ordering, i.e., $A^-$ consists
of  all $a\in A$ which act contractingly on the Lie algebra
${\n}$. Further let $\overline{A^-}$ be the closure of $A^-$ in
$G$, this is a manifold with boundary. Let $K_M$ be a maximal
compact subgroup of $M$. We may suppose that $K_M=M\cap K$ and
that $K_M$ contains $B$. Fix an irreducible unitary
representation $(\tau ,V_\tau)$ of $K_M$.

Let $\CE_P(\Ga)$ denote the set of all conjugacy classes $[\ga]$
in $\Ga$ such that $\ga$ is in $G$ conjugate to an element $a_\ga
b_\ga$ of $A^- B$.

Take a class $[\ga]$ in $\CE_P(\Ga)$. Then there is a conjugate
$H_\ga =A_\ga B_\ga$ of $H$ that contains $\ga$. Then the
centralizer $\Ga_{\ga}$ projects to a lattice $\Ga_{A,\ga}$in the
split part $A_\ga$. Let $\la_\ga$ be the covolume of this lattice.

Let $(\omega,V_\omega)$ be a finite dimensional unitary
representation of $\Ga$.

For any module $V$ of the Lie algebra ${\n}$ let $H^q({\n},V)$,
$q=0,\dots ,\dim\n$ denote the {\it Lie algebra cohomology}
 \cite{BorWall}. If $\pi\in\hat{G}$
then $H^q({\n},\pi_K)$ is an admissible $(\a\oplus\m,K_M)$-module
of finite length \cite{HeSch}.

For $\mu\in\a^*\cong\Hom(A,\C^\times)$ and $j\in\N$ let $C^{\mu,j}(A^-)$ denote the space of all
functions on $A$ which
\begin{itemize}
\item
are $j$-times continuously differentiable on $A$,
\item
are zero outside $A^-$,
\item
satisfy $|D\ph|\ll  |a^\mu|$ for every invariant differential operator $D$
on $A$ of degree $\le j$. 
\end{itemize}

This space can be topologized with the seminorms
$$
N_D(\ph)\= \sup_{a\in A} |a^{-\mu}D\ph(a)|,
$$
$D\in U(\a)$, $\deg(D)\le j$. Since the space of operators $D$ as above is
finite dimensional, one can choose a basis $D_1,\dots,D_n$ and set
$$
\norm{\ph}\= N_{D_1}(\ph)+\cdots +N_{D_n}(\ph).
$$
The topology of $C^{\mu,j}(A^-)$ is given by this norm and thus
$C^{\mu,j}(A^-)$ is a Banach space.

\begin{theorem}\label{lefschetz}  (Lefschetz formula)
 Assume
that $M$ is orientation preserving or that $\tau$ lies in the
image of the restriction map $\res_{K_M}^M$. There is $\mu\in\a^*$ and $j\in\N$ such that for every $\ph\in C^{\mu,j}(A^-)$ the expression
$$
\sum_{\pi \in \hat{G}} N_{\Ga,\omega} (\pi) \sum_{p,q}(-1)^{p+q}
 \int_{A^-} \ph (a)
\tr\(a|\(H^q(\n,\pi_K)\otimes\wedge^p\p_M \otimes
V_{\breve{\tau}}\)^{K_M}\)da,
$$
henceforth referred to as the global side, equals
$$
(-1)^{\dim N}\sum_{[\ga]\in {\cal E}_P(\Ga)}\la_\ga \
\chi_r(\Ga_\ga)\ \tr\omega(\ga)\ \frac{\ph(a_\ga)\tr
\tau(b_\ga)}{\det(1- a_\ga b_\ga |{\n})},
$$
called the local side, where $r=\dim A$.
Either side defines a continuous linear functional on the Banach space $C^{\mu,j}(A^-)$.
\end{theorem}

The proof will be given in section \ref{prf-lef}

We will give a reformulation of the theorem and for this
we need the following notation. Let $V$ be a complex
vector space on which $A$ acts linearly. For each
$\la\in\a^*$ let $V^\la$ denote the generalized
$\la$-eigenspace in $V$, i.e., $V^\la$ consists of all
$v\in V$ such that there is $n\in \N$ with
$$
(a-a^{\la})^n v\= 0
$$
for every $a\in A$. The theorem above implies  the
following.

\begin{theorem}\label{4.2}
Assume that $M$ is orientation preserving or that $\tau$
lies in the image of the restriction map $\res_{K_M}^M$.
Then we have the following identity of distributions on
$A^-$.
$$
\sum_{\pi\in\hat G} N_{\Ga,\omega}(\pi)\sum_{\la\in\a^*}m_\la(\pi)\, (\cdot)^\la\=
\sum_{[\ga]\in\CE_P(\Ga)} c_\ga \delta_{a_\ga}.
$$
Here $(\cdot)^\la$ means the function $a\mapsto a^\la$ and
 $m_\la(\pi)$ equals
$$
\sum_{p,q}(-1)^{p+q+\dim N}
\dim\left(H^q(\n,\pi_K)^\la\otimes\wedge^p\p_M \otimes
V_{\breve{\tau}}\right)^{K_M}.
$$
The sum indeed is finite for each $\la\in\a^*$. Further,
for $[\ga]\in {\cal E}_P(\Ga)$ we set
$$
c_\ga\=\la_\ga \ \chi_r(\Ga_\ga)\ \tr\omega(\ga)\
\frac{\tr \tau(b_\ga)}{\det(1- a_\ga b_\ga |{\n})}.
$$
\end{theorem}

\begin{corollary} \label{ruelle-lef}
Assume
that $M$ is orientation preserving or that $\tau$ lies in the
image of the restriction map $\res_{K_M}^M$. There is $\mu\in\a^*$ and $j\in\N$ such that for every $\ph\in C^{\mu,j}(A^-)$ the expression
$$
\sum_{\pi \in \hat{G}} N_{\Ga,\omega} (\pi)
\sum_{p,q,r}(-1)^{p+q+r+\dim N}\times\hspace{150pt}
$$ $$
\hspace{50pt} \int_{A^-} \ph (a)
\tr(a|(H^q(\n,\pi_K)\otimes\wedge^p\p_M \otimes\wedge^r\n\otimes
V_{\breve{\tau}})^{K_M})da,
$$
equals
$$
\sum_{[\ga]\in {\cal E}_P(\Ga)}\la_\ga \ \chi_r(\Ga_\ga)\
\tr\omega(\ga)\ \ph(a_\ga)\tr \tau(b_\ga).
$$
\end{corollary}

\prf Let $\wedge^r\n=\bigoplus_{j\in I_r} V_j$ be the
decomposition of the adjoint action of the group $M$ into
irreducible representations. On $V_j$ the torus $A$ acts
by a character $\la_j$. We apply the theorem to $\tau$
replaced by $\tau\oplus V_j|_{K_M}$ and $\ph(a)$ replaced
by $\ph(a)\la_j(a)$. We sum these over $j$ and take the
alternating sum with respect to $r$. On the local side we
apply the identity
$$
\sum_{r=0}^{\dim N}(-1)^r\tr (a_\ga b_\ga |\wedge^r\n) =
\det(1-a_\ga b_\ga|\n)
$$
to get exactly the local side of the corollary. On the global side
we need to recall that the $K_M$-module $\wedge^r\n$ is self dual.
The corollary follows. \qed

\subsection{Proof of the Lefschetz formula}\label{prf-lef}
The notations are as in section \ref{4.1}. Let $G$ act on itself
by conjugation, write $g.x = gxg^{-1}$, write $G.x$ for the orbit,
so $G.x = \{ gxg^{-1} | g\in G \}$ as well as $G.S = \{ gsg^{-1} |
s\in S , g\in G \}$ for any subset $S$ of $G$. We are going to
consider functions that are supported in the set $G.(MA^-)$. By
Theorem \ref{exist-ep} there exists an Euler-Poincar\'e function
$f_\tau^M\in C_c^\infty(M)$ to the representation $\tau$.

For a finite dimensional complex vector space $V$ and $ T\in \GL(V)$ let $E(T)$ be the set of eigenvalues of $T$.
Let $\la_{min}(T):= \min\{ |\la| : \la\in E(T)\}$ and  $\la_{max}(T):=
\max\{ |\la| : \la\in E(T)\}$.
We are particularly interested in the adjoint action of $G$ on its Lie algebra $\g$.
So for $g\in G$ and $V$ a $g$-invariant subspace of $\g$ we write $g|V$ for the induced element of $\GL(V)$.

For $am\in AM$ define
$$
\la(am) := \frac{\la_{min}(a|\bar{\n})}{\la_{max}(m|\g)}.
$$
Note that $\la_{max}(m|\g)$ is always $\ge 1$ and that
$\la_{max}(m|\g)\la_{min}(m|\g)=1$. 

We will consider the set
$$
(AM)^\sim \ :=\ \{ am\in AM | \la(am)>1 \}.
$$
Let $M_{ell}$ denote the set of elliptic elements in $M$.

\begin{lemma} \label{MA}
The set $(AM)^\sim$ has the following properties:

\begin{itemize}
\item[1.]
$A^-M_{ell}\subset (AM)^{\sim}$
\item[2.]
$am\in (AM)^\sim \Rightarrow a\in A^-$
\item[3.]
$am, a'm' \in (AM)^\sim, g\in G\ {\rm with}\ a'm'=gamg^{-1}
\Rightarrow a=a', g\in AM$.
\end{itemize}
\end{lemma}

\prf The first two are immediate. For the third let $am, a'm' \in
(AM)^\sim$ and $g\in G$ with $a'm'=gamg^{-1}$. Observe that by the
definition of $(AM)^\sim$ we have
\begin{eqnarray*}
\la_{min}(am | \bar{\n}) &\ge& \la_{min}(a| \bar{\n})\la_{min}(m |
\g)\\
 &>& \la_{max}(m|\g)^2\la_{min}(m|\g)\\
 &=& \la_{max}(m|\g)\\
 &\ge & \la_{max}(m | \a +\m +\n)\\ &\ge&\la_{max}(am | \a +\m +\n)
\end{eqnarray*}
that is, any eigenvalue of $am$ on $\bar{\n}$ is strictly bigger
than any eigenvalue on $\a +\m +\n$. Since $\g = \a +\m +\n
+\bar{\n}$ and the same holds for $a'm'$, which has the same
eigenvalues as $am$, we infer that $\Ad(g)\bar{\n} =\bar{\n}$. So
$g$ lies in the normalizer of $\bar{\n}$, which is
$\bar{P}=MA\bar{N} =\bar{N}AM$. Now suppose $g=nm_1a_1$ and
$\hat{m} =m_1mm_1^{-1}$ then
$$
gamg^{-1} \= na\hat{m}n^{-1} \= a\hat{m}\
(a\hat{m})^{-1}n(a\hat{m})\ n^{-1}.
$$
Since this lies in $AM$ we have $(a\hat{m})^{-1}n(a\hat{m}) =n$
which since $am\in (AM)^\sim$ implies $n=1$. The lemma is proven.
\qed

Let $C\subset M$ be a compact subset.
In our application, $C$ will be the support of the Euler-Poincar\'e function $f_\tau^M$.
Let $\{\al_1,\dots,\al_r\}$ be the set of simple roots in $\phi^-(\g,\a)$.
Then $\al_1,\dots,\al_r$ is a basis of $\a^*$ and $\a_\R^-$ is the set of all $X\in\a$ with $\al_j(X)>0$ for $j=1,\dots,r$.
For $a\in A$ we write $a_j=\al_j(\log a)$ and thus we get global co-ordinates on $A$ such that $a\in A^-\Leftrightarrow a_j>0\ \forall j$.
For $T>0$ set $A_T=\{ a\in A : a_j\ge T\ \forall j\}$.
Then there exists $T>0$ such that the set closed $A_T\cdot C$ is contained in $(AM)^\sim$.

The boundary $S$ of $(AM)^\sim$ in $AM$ decomposes into two disjoint subsets $S=S_1\cup S_2$, where
\begin{eqnarray*}
S_1 &=& \{ am\in S : a\ne 1\}\\
S_2 &=& \{ am\in S : a=1\}.
\end{eqnarray*}

We want to construct a smooth function $\chi\colon (AM)^\sim\to [0,1]$ such that
\begin{itemize}
\item $\chi$ vanishes to infinte order at every point of $S_1$.
\item $\chi$ is invariant under conjugation by elements of $m$.
\item $\chi(am)=1$ if $m$ is elliptic.
\item If $a_j\ge T$ for some $j$ and $m\in C$, then 
$$
\frac\partial {\partial a_j}\chi(am)\= 0.
$$
\end{itemize}
Note that the last condition implies that $\chi$ is constant on $A_T\cdot C$.

In order to construct conjugation invariant functions one considers the geometric quotient $M/\conj$ which is an affine variety as $M$ is reductive.
Note that on the complex valued points the map $M_\C\to M/\conj(\C)$ is open as a consequence of the Kempf-Ness Theorem \cite{KempfNess}.
Embed $M/\conj$ into affine space $\A^n$, then $M/\conj(\C)\hookrightarrow \A^n(\C)\cong\C^n$.
We use the isomorphism of $A$ with $\a_\R$, so we embed $A$ into $\C^r$.
Thus we get a map
$$
\al \colon AM\to AM/\conj = a\times (M/\conj)\hookrightarrow C^{n+r}.
$$
Let $(A_\C M_\C)^\sim$ be the set of all $am\in A_\C M_\C$ with $\la(am)>1$.
Then $(A_\C M_\C)^\sim$ is a complex neighbourhood of $(AM)^\sim$ and there is an open subset $U$ of $\C^{n+r}$ such that $((A_\C M_\C)^\sim=\al^{-1}(U\cap\al(A_\C M_\C))$.
It follows that for each compact subset $C$ of $\C^n$ there exists $T>0$ such that for every $z\in\C^r$ with $\Re(z_j)\ge T\ \forall j$ one has $C\times \{ z\}\subset U$.
The task of constructing $\chi$ now boils down to constructing a function on $U$ with the indicated properties which is easily established.

Extend $\chi$ from $(AM)^\sim$ to all of $AM$ by setting
$$
\chi(am)\= 0,\quad {\rm if}\ am\notin (AM)^\sim.
$$
Fix a smooth function $\eta$ on $N$ which has compact support, is
positive, invariant under $K_M$ and satisfies $\int_N\eta(n) dn
=1$. 
Given these data let
$f = f_{\eta ,\tau ,\ph} : G\ra \C$ be defined by
$$
f (kn ma (kn)^{-1}) := \eta (n) f_\tau^M(m)
\frac{{\ph}(a)\chi(am)}{\det(1-(ma)|\n)},
$$
for $k\in K, n\in N, m\in M, a\in{A}$. Further $f(x)=0$
if $x$ is not in $G.(AM)$.
Note that indeed, $f$ is supported in the closure of $G.(AM)^\sim$.

\begin{lemma}
The function $\frac{\ph(a)\chi(am)}{\det(1-am |\n)}$ is $j$-times continuously differentiable on $(AM)^\sim$ and vanishes on the boundary $\partial (AM)^\sim$ in $AM$ to order at least
$j-\dim\n$.
\end{lemma}

\prf
Let $a_0m_0$ be a boundary point of $(AM)^\sim$.
If $a_0\in A^-$, then $\chi$ vanishes at $a_0m_0$ to infinite order, and so does  $\frac{\ph(a)\chi(am)}{\det(1-am |\n)}$.
If $a_0$ lies on the boundary of $A^-$, then the vanishing order is determined by $\frac{\ph(a)}{\det(1-am|\n)}$.
\qed

\begin{lemma} \label{welldef}
The function $f$ is well defined and for given $N\in\N$ there are $\mu,j$ such that the map $C^{\mu,j}(A^-)\to L_{2N}^1(G)$; $\ph\mapsto f_{\eta,\tau,\ph}$ is continuous.
\end{lemma}

\prf By the decomposition $G=KP=KNMA$ every element $x\in
G.(AM)^\sim$ can be written in the form $kn ma (kn)^{-1}$. Now
suppose two such representations coincide, that is
$$
kn ma (kn)^{-1}\= k'n' m'a' (k'n')^{-1}
$$
then by Lemma \ref{MA} we get $(n')^{-1} (k')^{-1}kn\in MA$, or
$(k')^{-1}k\in n'MAn^{-1}\subset MAN$, hence $(k')^{-1}k\in K\cap
MAN=K\cap M=K_M$. Write $(k')^{-1}k=k_M$ and $n''=k_Mnk_M^{-1}$,
then it follows
$$
n'' k_Mmk_M^{-1} a (n'')^{-1}\= n'm'a'(n')^{-1}.
$$
Again by Lemma \ref{MA} we conclude $(n')^{-1}n''\in MA$, hence
$n'=n''$ and so
$$
k_Mmk_M^{-1} a \= m'a',
$$
which implies the well-definedness of $f$.

Let $N\in\N$. We will show that for $\Re(\mu),j$ and $l$ sufficiently large  the function $f$
lies in
$L^1_{2 N,1}(G)$.

Let the group $K_M=K\cap M$ act from the right on $K\times N\times M\times A$ by
$$
(k,n,m,a)k_M\= (kk_M,k_M^{-1}nk_M,k_M^{-1}mk_M,a).
$$
Let $(K\times N\times M\times A)/K_M$ denote the quotient, then the projection
$$
K\times N\times M\times A\ \to\ (K\times N\times M\times A)/K_M
$$
is a principal $K_M$-fibre bundle.

Consider the map
$$
\begin{array}{c c c c}
F\ :& (K\times N\times M\times A)/K_M &\ra& G\\
{}& [k,n,m,a]&\mapsto&k n a m(k n)^{-1}.
\end{array}
$$
Then $f$ is a $j-1-\dim(\n)$-times continuously differentiable
function on $Z=(K\times N\times M\times A)/K_M$ which factors over $F$.
Now set $D=(K\times N\times (AM)^\sim)/K_M\subset Z$.
Then $D$ is an open subset of $Z$.
Set $S=Z\setminus D$.
The first part of this proof shows that $F$ is a diffeomorphism on $D$ with open image.

Thus to compute the order of differentiability of $f$ we can apply Proposition \ref{smoothness}.
To compute the order of differentiability of $f$ as a function on $G$ we
have to take into count the zeroes of the differential of $F$. So
we compute the differential $TF$ of $F$. Let at first $X\in \k$, then
$$
TF(X) f(k n a m(k n)^{-1}) = \frac{d}{d t}|_{t=0} f(k
\exp(t X)n a m n^{-1}\exp(-t X)k^{-1}),
$$
which implies the equality
$$
TF(X)_x = (\Ad(k)(\Ad(n(am)^{-1}n^{-1})-1)X)_x,
$$
when $x$
equals $k n a m(k n)^{-1}$. 
Note that for $X\in\k$, unless $X\in\k_M$, we have $(n\Ad(am)n^{-1}-1)(X)\ne 0$.
Similarly for $X\in\n$ we get that
$$
TF(X)_x = (\Ad(k n)(\Ad((am)^{-1})-1)X)_x
$$
and for $X\in
\a\oplus \m$ we finally have $TF(X)_x = (\Ad(k n)X)_x$.
From this it becomes clear that $F$, regular on $K\times
N\times M \times A^-$, may on the boundary have vanishing
differential of order at most $\dim(\n)+\dim(\k)$. 
In order to apply Proposition \ref{smoothness} we next need to determine the vanishing order of $f\circ F$ at $S$.

Applying Proposition \ref{smoothness} we get
that $f$ is $\left[\frac{j-\dim\n-2}{\dim\n+\dim\k}\right]$-times continuously differentiable on $G$. So we assume
$j\ge 2\dim(\n)+\dim\k+2$ from now on. 
 We have to show that 
$D f\in L^1(G)$ for any $D\in U(\g)$ of degree $\le 2 N$.
The same computation will also show the boundedness of $Df$ for $\deg(D)\le 1$.
 For this we recall the map $F$ and our
computation of its differential. 
Let $\q\subset \k$ be a
complementary space to $\k_M$. 
On the regular set $TF$ is bijective. 
Fix $x=k n a m(k n)^{-1}$ in the
regular set and let $TF_{x}^{-1}$ denote the inverse map of $TF_{x}$
which maps to $\q\oplus\n\oplus\a\oplus\m$. Introducing norms on
the Lie algebras we get an operator norm for $TF_{x}^{-1}$ and the
above calculations show that $\norm{TF_{x}^{-1}}\le P(am)$, where $P$ is
a class function on $AM$, which, restricted to any Cartan $H=AB$
of $AM$ is a linear combination of quasi-characters. Supposing
$j$ and $\Re(\mu)$ large enough we get for $D\in U(\g)$ with
$\deg(D)\le 2 N$:
$$
|D f(k n a m(k n)^{-1})| \le \sum_{D_1} P_{D_1}(am) |D_1 f(k,n,a,m)|,
$$
where the sum runs over a finite set of $D_1\in
U(\k\oplus\n\oplus\a\oplus\m)$ of degree $\le 2 N$ and $P_{D_1}$
is a function of the type of $P$. On the right hand side we have
considered $f$ as a function on $K\times N\times A\times M$.
This discussion uses the facts that $K$ is compact, $N$ is unipotent, and $\det\left(\Ad(n(am)^{-1}n^{-1})-1\right)=\det\left(\Ad(am)^{-1}-1\right)$.
Finiteness of the sum in the inequality above follows from the Poincar\'e-Birkhoff-Witt Theorem.
Integrating, it becomes clear that for $\mu$ and $j$ sufficiently large the map $\ph\mapsto f$ indeed is continuous.
\qed

We will plug $f$ into the trace formula. For the geometric side
let $\ga \in \Ga$. We have to calculate the orbital integral:
$$
\CO_\ga (f) = \int_{G_\ga \bs G} f(x^{-1}\ga x) dx.
$$
by the definition of $f$ it follows that $\CO_\ga(f)=0$ if
$\ga\notin G.(AM)^\sim$. It remains to compute $\CO_{am}(f)$ for
$am\in(AM)^\sim$. Again by the definition of $f$ it follows
$$
\CO_{am}(f) \=
\CO_m^M(f_\tau^M)\frac{\tilde{\ph}(am)}{\det(1-ma|\n)},
$$
where $\CO_m^M$ denotes the orbital integral in the group $M$.

Since only elliptic elements have nonvanishing orbital integrals
at $f_\tau^M$ it follows that only those conjugacy classes $[\ga]$
contribute for which  $\ga$ is in $G$ conjugate to $a_\ga b_\ga\in
A^-B$. Recall that Theorem \ref{3.3} says
$$
\vol(\Ga_\ga\bs G_\ga)\= \chi_r(\Ga_\ga) \la_\ga
 \frac{c_{G_\ga}}
     {|W_{\ga,\C}|\prod_{\alpha\in\Phi_\ga^+}(\rho_\ga,\alpha)}.
$$
By Theorem \ref{orbitalint} we on the other hand get
$$
\CO_\ga(f) \=
\frac{|W_{\ga,\C}|\prod_{\alpha\in\Phi_\ga^+}(\rho_\ga,\alpha)}
     {c_{G_\ga}}
     \tr\tau(b_\ga)\frac{\ph(a_\ga)}
                        {\det(1-a_\ga t_\ga|\n)},
$$
so that
$$
\vol(\Ga_\ga\bs G_\ga)\, \CO_\ga(f)\= \chi_r(\Ga_\ga) \la_\ga
\tr\tau(b_\ga)\frac{\ph(a_\ga)}
                        {\det(1-a_\ga b_\ga|\n)}.
$$
It follows that the geometric side of the trace formula coincides
with the geometric side of the Lefschetz formula.

Now for the spectral side let $\pi\in\hat{G}$. We want to compute
$\tr\pi(f)$. Let $\Theta_\pi^G$ be the locally integrable
conjugation invariant function  on $G$ such that
$$
\tr\pi(f)\= \int_G f(x) \Theta_\pi^G(x) dx.
$$
To evaluate $\tr \pi(f)$ we will employ the Hecht-Schmid
character formula \cite{HeSch}. For this let
$$
(AM)^- = {\rm interior\ in\ } MA {\rm \ of\ the\ set}
$$ $$
\left\{ g\in MA | \det (1-ga | \n) \geq 0\ {\rm for\ all\ } \a\in
A^- \right\}.
$$
Note that $(AM)^\sim$ is a subset of $(AM)^-$.
The main result of \cite{HeSch} is that for $ma\in (AM)^- \cap
G^{reg}$, the regular set, we have
$$
\Theta_\pi^G(am) = \frac{\sum_{p=0}^{\dim \n}(-1)^p \Theta_{H_p(\n
,\pi_K)}^{MA}(am)}{\det (1-am |\n)}.
$$
Let $h$ be supported on $G.(AM)^-$, then the Weyl integration
formula implies that
$$
\int_G f(x) dx = \int_{G/MA} \int_{MA^-} h(gmag^{-1})
|\det(1-ma|\n\oplus\bar{\n})\,dadmdg.
$$
So that for $\pi\in\hat{G}$:
\begin{eqnarray*}
\tr \pi(f) &=& \int_G \Theta_\pi^G (x) f (x) dx\\
    &=& \int_{MA^-} \Theta_\pi^G (ma) f_\tau^M(m)\tilde{\ph}(am)
|\det(1-ma|\bar{\n})|\,dadm\\
    &=&  \int_{MA^-}f_\tau^M(m)\sum_{p=0}^{\dim
    N}(-1)^p\Theta_{H_p(\n,\pi_K)}^{AM}(am)\\
    &{}&\hspace{20pt}\times
    \frac{|\det(1-am|\bar{\n})|}
         {\det(1-am|\n)}
    \tilde{\ph}(am)\ dadm.
\end{eqnarray*}
Now we find that
\begin{eqnarray*}
|\det(1-am|\bar{\n})| &=& (-1)^{\dim N} \det(1-am|\bar{\n})\\
    &=& (-1)^{\dim N} a^{-2\rho_P}\det(a^{-1}-m|\bar{\n})\\
    &=& (-1)^{\dim N} a^{-2\rho_P}\det((am)^{-1}-1|\bar{\n})\\
    &=& a^{-2\rho_P}\det(1-(am)^{-1}|\bar{\n})\\
    &=& a^{-2\rho_P}\det(1-am|{\n})
\end{eqnarray*}
so that
$$
\tr\pi(f)\= \int_{MA^-}f_\tau^M(m)\sum_{p=0}^{\dim
N}(-1)^p\Theta_{H_p(\n,\pi_K)}^{AM}(am)a^{-2\rho_P}\tilde{\ph}(am)dadm.
$$
We have an isomorphism of $(\a\oplus\m,K_M)$-modules \cite{HeSch}
$$
H_p(\n,\pi_K)\cong H^{\dim N-p}(\n,\pi_K)\otimes \wedge^{top}\n.
$$
This implies
$$
\sum_{p=0}^{\dim N}(-1)^p\Theta_{H_p(\n,\pi_K)}^{AM}(am)
a^{-2\rho_P}\= (-1)^{\dim N} \sum_{p=0}^{\dim
N}(-1)^p\Theta_{H^p(\n,\pi_K)}^{AM}(am).
$$
And so
$$
\tr\pi(f)\= \int_{MA^-}f_\tau^M(m) \sum_{p=0}^{\dim N}(-1)^{p+\dim
N} \Theta_{H^p(\n,\pi_K)}^{AM}(am)\tilde{\ph}(am)\ dadm.
$$
Let $B=H_1,\dots,H_n$ be the conjugacy classes of Cartan subgroups
in $M$. By the Weyl integration formula the integral over $M$ is a
sum of expressions of the form
$$
\int_{H_j}\int_{M/H_j}
f_\tau^M(x^{-1}hx)\Theta_{H^*(\n,\pi_K)}^{MA}(x^{-1}hax)
\tilde{\ph}(x^{-1}hax)dx |\det(1-h|\m/\h_j)|dh
$$ $$
=\int_{H_j}\int_{M/H_j}
f_\tau^M(x^{-1}hx)\Theta_{H^*(\n,\pi_K)}^{MA}(ha)
\tilde{\ph}(ha)dx |\det(1-h|\m/\h_j)|dh
$$ $$
=\int_{H_j}\CO_h^M(f_\tau^M)\Theta_{H^*(\n,\pi_K)}^{MA}(ha)
\tilde{\ph}(ha) |\det(1-h|\m/\h_j)|dh,
$$
where we have used the conjugacy invariance of
$\Theta_{H^*(\n,\pi_K)}^{MA}$ and $\tilde{\ph}$. The orbital
integrals $\CO_h^M(f_\tau^M)$ are nonvanishing only for $h$
elliptic, so only the summand with $H_j=H_1=B$ survives. In this
term we may replace $\tilde{\ph}(ha)$ by $\ph(a)$ so that we get
$$
\tr\pi(f)\= \int_{MA^-}f_\tau^M(m)\sum_{p=0}^{\dim N}(-1)^{p+\dim
N} \Theta_{H^p(\n,\pi_K)}^{MA}(am) \ph(a)da dm
$$ $$
= \sum_{p,q\ge 0} (-1)^{p+q+\dim N}\int_{A^-} \ph(a)
\tr(a|(H^q(\n,\pi_K)\otimes\wedge^p\p_M\otimes
V_{\breve{\tau}})^{K_M}) da.
$$
The theorem follows.
\qed

\end{document}